\documentclass{article}
\usepackage{amsmath,amstext,amsthm,amsfonts}
\usepackage[dvips]{graphicx}
\usepackage{epic,eepic}

\theoremstyle{plain}
\newtheorem{theorem}{Theorem}[section]
\newtheorem{proposition}[theorem]{Proposition}
\newtheorem{lemma}[theorem]{Lemma}
\newtheorem{corollary}[theorem]{Corollary}
\newtheorem{maintheorem}{Theorem}

\theoremstyle{definition}
\newtheorem{remark}[theorem]{Remark}

\newtheorem{claim}[theorem]{Claim}
\newtheorem{definition}[theorem]{Definition}

\newcommand{\field}[1]{\mathbb{#1}}
\newcommand{\real}{\field{R}}

\renewcommand{\natural}{\field{N}}

\newcommand{\al} {\alpha}       
        
\newcommand{\ga} {\gamma}    
\newcommand{\de} {\delta}       \newcommand{\De}{\Delta}
\newcommand{\ep} {\epsilon}
\newcommand{\vep}{\varepsilon}

\newcommand{\ka} {\kappa}

\newcommand{\SB}{{\cal B}}
\newcommand{\SC}{{\cal C}}

\newcommand{\SE}{{\cal E}}
\newcommand{\SF}{{\cal F}}

\newcommand{\SL}{{\cal L}}

\newcommand{\SN}{{\cal N}}

\newcommand{\SP}{{\cal P}}
\newcommand{\SQ}{{\cal Q}}
\newcommand{\SR}{{\cal R}}

\begin{document}

\title{Fast decay of correlations of equilibrium states of 
open classes of non-uniformly expanding maps and potentials}
\author{Alexander Arbieto and Carlos Matheus}
\date{March 27, 2006}

\maketitle

\begin{abstract}
We study the existence, uniqueness and rate of decay of correlation of 
equilibrium measures associated to robust classes of non-uniformly expanding 
local diffeomorphisms and H\"older continuous potentials. The approach used in
this paper is the spectral analysis of the Ruelle-Perron-Frobenius transfer 
operator. More precisely, we combine the expanding features of the eigenmeasures
of the transfer operator with a Lasota-Yorke type inequality to prove the
existence of an unique equilibrium measure with fast decay of correlations.
\end{abstract}

\section{Introduction}

The thermodynamical formalism of uniformly expanding local diffeomorphisms
developed by Ruelle, Sinai, Bowen (among others) in the seventies is a major
example of applications of ideas of statistical mechanics into the ergodic
theory of dynamical systems. Indeed, the dictionary between one-dimensional
lattices and expanding systems (in particular, Gibbs distributions and 
equilibrium measures) established via Markov partitions permits to translate 
results in statistical mechanics into relevant theorems in the ergodic theory of
uniformly expanding systems. 

However, the extension of this theory beyond the uniformly expanding context
reveals fundamental difficulties, even if we are restricted to the non-uniformly
expanding realm (i.e., abundance of positive Lyapounov exponents). In fact, when
dealing with such problem, a major trouble is the ausence of generating finite 
Markov partitions: in general, being very optmistic, one can expect only 
Markov partitions with infinitely many symbols, leading us to consider the
thermodynamical formalism of gases with infinitely many states, a hard subject
not yet completely understood (see e.g.~\cite{BS} for a recent progress). 

On the other hand, some substantial advance was made by Alves, Bonatti and
Viana~\cite{ABV} concerning the existence and uniqueness of SRB measures for
some robust classes of non-uniformly expanding maps: they essentially assume 
that \emph{Lebesgue} almost every point has only non-zero Lyapounov exponents.
Since the physical measures should be absolutely continuous with respect to the
Lebesgue measure, this condition is really says something about the problem of
existence of SRB measures. But, it is not clear how to formulate the correct 
analog of this assumption in the case of equilibrium states with respect 
to arbitrary potentials (because, in general, equilibrium measures are not 
absolutely continuous with respect to Lebesgue). In this direction, some 
recent contributions have been given by several authors: see, e.g., the work of 
Bruin and Keller~\cite{BK} for interval maps, Denker and
Urbanski~\cite{DU},~\cite{U} for rational functions of the Riemann sphere,
Wang and Young~\cite{WY} for H\'enon-like maps; Buzzi, Maume-Deschamps, Sarig 
and Yuri~\cite{Bu},~\cite{BMD},~\cite{BS},~\cite{S},~\cite{Yu} for countable 
Markov shifts and piecewise expanding maps. Many of these papers deals with
dynamical systems with neutral periodic points, a relevant source of example of
non-hyperbolic systems. 

Nevertheless, we consider a different open class $\SC$ of non-uniformly 
expanding maps obtained through \emph{pitch-fork bifurcation of 
a periodic point of an expanding linear torus endomorphism}. Here, the maps are 
assumed to be uniformly
expanding outside a small neighborhood of the periodic point where the
pitch-fork bifurcation was made (but possibly of saddle type near the periodic
point) and the potential $\phi$ is supposed to have its
oscillation $\max\phi-\min\phi$ not too large. In this setting,
Oliveira~\cite{O} proved the existence of ergodic equilibrium states whose 
Lyapounov exponents are positive. Also, Arbieto, Matheus and Oliveira~\cite{AMO}
extended Oliveira's approach to obtain similar results for small random
pertubations of these maps. 

Recently, Oliveira and Viana~\cite{OV1} improved
the results of~\cite{O} when the potential $\phi$ is constant 
by showing that the equilibrium state with respect to 
this potential (i.e., the measure of maximal entropy) is unique. The idea is 
to use the \emph{Ruelle-Perron-Frobenius transfer operator} to construct a 
maximal measure. Then, the uniqueness result follows from the good 
spectral properties of Ruelle's operator. We would like to point out that
Oliveira and Viana annouced in the preprint~\cite{OV2} that this method can be
pushed a little bit more to extend the results of~\cite{OV1} for non-constant
potentials with low oscillation.

The present paper relies on the basic strategy of~\cite{OV2}, i.e., 
we analyze the spectral properties of the transfer operator obtaining our main
result: 

\begin{themaintheorem}``Given $f\in\SC$ and $\phi$ is an H\"older continuous
potential with $\max\phi-\min\phi$ sufficiently small, there exists an unique
equilibrium measure $\mu$ with respect to $(f,\phi)$. Also, the Lyapounov
exponents of $\mu$ are positive and $\mu$ has exponential decay of
correlations.''
\end{themaintheorem}

The organization of this paper is as follows. In section 2, we introduce the
open classes of non-uniformly expanding maps and H\"older continuous potentials
(with low variation) treated here. In section 3, we consider  the
Ruelle-Perron-Frobenius transfer operator $\SL_{\phi}$ and the eigenmeasures
$\nu$ of its dual operator $\SL_{\phi}^{\star}$. After that, we prove that $\nu$
is expanding (i.e., hyperbolic times are abundant at $\nu$-generic points) and
$\SL_{\phi}$ satisfies a Lasota-Yorke tye inequality. Next, we show that this
ensures the existence of an eigenfunction $h$ of $\SL_{\phi}$. As a consequence,
we can obtain the central result of the section 3: the existence of an unique
$f$-invariant probability $\mu$ which is absolutely continuous with respect to
$\nu$ (= an eigenmeasure of $\SL_{\phi}^{\star}$); moreover, $\mu$ has
exponential decay of correlations. In section 4, we prove that the $f$-invariant
probability $\mu$ is an equilibrium state of $(f,\phi)$ and any equilibrium
measure $\eta$ verifies that $\frac{1}{h}\eta$ is an eigenmeasure of
$\SL_{\phi}^{\star}$. In section 5, we complete the proof of our main result by
showing the uniqueness of equilibrium measure $\mu$. The idea is that the
expanding features of any eigenmeasure of $\SL_{\phi}^{\star}$ implies that they
are all equivalent so that $\frac{1}{h}\eta$ is equivalent to $\nu$ from the
results of section 4. Then, by ergodicity, it follows that $\eta=h\nu=\mu$, as
desired. Finally, in section 6, we include some ergodic properties of the
equilibrium states constructed here: spectral gap, stochastic stability and
large deviations. The proof of these results are simple modifications of their
analogous in the uniformly expanding context (based on the lemmas proved in this
article) and they were included only to illustrate the power of the
machinary developed in the previous sections. To make the exposition more
self-contained, we included also three appendices containing some facts used 
along the proof of the main result (but not showed in the text in order to do
not interrupt the argument).

To close the introduction, we would like to say that this article is only starts
a program to understand the ergodic theory of non hyperbolic systems. Indeed,
some questions motivated by our theorems are: 
\begin{itemize}
\item Is it possible to extend our existence, uniqueness and decay
of correlations results for more general potentials? 
\item What about analogous theorem for robust classes of \emph{non-uniformly
hyperbolic diffeomorphisms} (e.g., pitch-fork bifurcation of the expanding
direction of a periodic saddle-point of an Anosov diffeomorphism) and/or 
\emph{non-uniformly expanding local diffeomorphims with critical points} (e.g.,
Viana maps)?
\end{itemize}  
 
\section{Preliminaries}

Let $f:M\to M$ be a continuous transformation on a compact space $M$ and
$\phi:M\to\real$ be a continuous function. An $f$-invariant measure is an
\emph{equilibrium state} of $f$ for the potential $\phi$ if it maximizes the
functional 
\begin{equation*}
\eta\mapsto h_{\eta} (f) + \int\phi d\eta,
\end{equation*}
among all $f$-invariant probabilities $\eta$.

\subsection{The class of transformations and potentials}

We consider $f:M\to M$ is a $C^1$ local diffeomorphism on a compact 
boundaryless manifold $M$ such that

\bigskip

(H1) There exists connected closed sets (with finite inner diameter) 
$$\SR = \{R_1,\dots ,R_q,R_{q+1},\dots ,R_{p+q}\}$$
whose interiors are pairwise disjoint such that $\text{diam}(\SR)<\vep$ for some
$\vep>0$, $\bigcup\limits_{i} R_i = M$, 
$f$ is injective on each
$R_i$ and, for some constants $\sigma_1 < 1$ and $\delta_0 > 0$,
\begin{itemize}
\item $f$ is expanding at every 
$x\in R_{q+1}\cup\dots\cup R_{p+q}$: $\|Df(x)^{-1}\|\leq
\sigma_1^{-1}$;
\item $f$ is never too contracting: $\|Df(x)^{-1}\|\leq 1+\delta_0$ for every
$x\in M$;
\item $\SR$ is a transitive Markov partition: the image $f(R_i)$ of
every atom is the union of some atoms $R_j$ and there exists $N$ such that
$f^N(R_i)=M$, $i=1\dots,p+q$;
\end{itemize}

It is worth to point out that the partition $\SR$ \emph{is not required to be
generating}, i.e., the diameters of the cylinders of lenght $n$ do not need to
decrease to zero when $n\to\infty$.

We denote by $k$ the number of pre-images of any point $x$ under $f$. This
number does not depend on $x$ because $f$ is a local diffeomorphism. We assume
that $k>q$ (i.e., every point has some pre-image in a good rectangle $R_i$,
$i=q+1\dots,p+q$). Moreover, the potential $\phi: M\to\real$ is supposed to
verify 

\bigskip

(H2) $\phi$ is $\al$-H\"older-continuous and its variation 
$\max\phi - \min\phi$ is sufficiently small, i.e., 
$$
\max\phi - \min\phi < \ep_0 (f).
$$

\begin{remark}By the definition of equilibrium states, it is not difficult to
see that any equilibrium state $\eta$ of $\phi$ is also an equilibrium state of
the potential $\phi+c$, for every constant $c\in\real$. So, up to consider 
$\phi - \inf\phi$ instead of $\phi$, we can assume that $\inf\phi=0$ and, 
in particular, $\phi\geq 0$.  
\end{remark} 

\subsection{Statement of the main theorems} 

The main theorems of this paper concerns the existence, uniqueness, fast mixing
(exponential decay of correlations) and the oscillation of the Birkhoff sums of
H\"older observables around the expected value (central limit theorem).
Unfortunately, it turns out that our techniques only can be applied if some
extra control of the contraction of the derivative of $f$ is assumed. 

In this article, we deal with two different candidates to be the extra 
condition, both of them ensuring the necessary control of the contration for 
the application of our method.

The first candidate is:

\bigskip

(H3) $\max\limits_{x\in M}\log\|\Lambda^{d-1} Df(x)\| < \log k$, where $d$ is 
the dimension of $M$ and $\Lambda^k$ is the $k$-th exterior power of $Df$.

\bigskip

Morally speaking, (H3) means that the $(d-1)$-dimensional volume is not expanded
too much.

Our first main result is related with the ergodic theory of equilibrium states 
of the class of transformations verifying (H1), (H2) and (H3):

\begin{maintheorem}\label{t.A}Under the hypothesis (H1), (H2) and (H3) above, 
if $\delta_0$ is sufficiently small, then $f$ has an unique equilibrium state
$\mu$ for the potential $\phi$. Moreover, $\mu$ has exponential decay of
correlations for H\"older observables:
there exists some constants $0<\tau<1$ and $K>0$ such that,
for all $u\in L^1(\nu), v\in C^{\alpha}(M)$ and $n\geq 0$,
\begin{equation*}
\left|\int_M (u\circ f^n)v d\mu - \int_M u d\mu\int_M v d\mu\right|\leq 
K(u,v)\cdot\tau^n,
\end{equation*} 
and the central limit theorem property:
for all $u$ H\"older continuous,
the random variable 
\begin{equation*}
\frac{1}{\sqrt{n}}\sum\limits_{j=0}^{n-1}(u\circ f^j - \int_M u d\mu) 
\end{equation*}
converges to a Gaussian law. More precisely, let $u$ be a $\alpha$-H\"older 
continuous function and  
\begin{equation*}
\sigma^2:=\int v^2 d\mu + 2\sum\limits_{j=1}^{\infty} v\cdot (v\circ f^j) d\mu,
\quad \text{ where } \quad v=u-\int u d\mu.
\end{equation*}
Then $\sigma<\infty$ and $\sigma=0$ iff $u=\varphi\circ f - \varphi$ for some
$\varphi\in L^1(\mu)$. If $\sigma>0$, then, for every interval $A\subset\real$, 
\begin{equation*}
\mu\left(x\in M: \frac{1}{\sqrt{n}}\sum\limits_{j=0}^{n-1}
\left(u(f^j(x)-\int u d\mu)\right)\in A\right)\to \frac{1}{\sigma\sqrt{2\pi}}
\int_A e^{-\frac{t^2}{2\sigma^2}} dt,
\end{equation*}
as $n\to\infty$.
\end{maintheorem} 

The second candidate is:

\bigskip

(H4) $|\det Df(x)|\geq\sigma_2>q$ for all $x\in M$ (i.e., $f$ is
volume-expanding everywhere) and there exists a set 
$W\subset R_1\cup\dots\cup R_q$ containing 
$V:=\{x\in M: \ \|Df(x)^{-1}\|>\sigma_1^{-1}\}$ such that 
$$ M_1 > m_2 \quad \text{ and } \quad m_2 - m_1 < \beta,$$
where $m_1$, $m_2$ are the infimum and the supremum of $\log\|\det Df\|$ on
$V$, resp. and $M_1$, $M_2$ are the infimum and the supremum of $\log\|\det
Df\|$ on $M-W$, resp.

\bigskip

Informally, (H4) says that $f$ is volume-expanding everywhere, i.e., $f$ can
not contract all the directions at the same time, so that the basins of 
sinks are not allowed as a region where $f$ can eventually contract. Also, the 
expansion of volume near the region of $Df$-contraction is close to be constant
and, furthermore, the volume is expanded more strongly outside the region of
contraction than inside of it.  

Our second main result is related with the ergodic theory of equilibrium states 
of the class of transformations verifying (H1), (H2) and (H4): 

\begin{maintheorem}\label{t.B}Under the hypothesis (H1), (H2) and (H4) above, 
if $\beta$ and $\delta_0$ are sufficiently small, then $f$ has an unique 
equilibrium state $\mu$ for the potential $\phi$. Moreover, $\mu$ has 
exponential decay of correlations for H\"older observables: 
there exists some constants $0<\tau<1$ and $K>0$ such that,
for all $u\in L^1(\nu), v\in C^{\alpha}(M)$ and $n\geq 0$,
\begin{equation*}
\left|\int_M (u\circ f^n)v d\mu - \int_M u d\mu\int_M v d\mu\right|\leq 
K(u,v)\cdot\tau^n,
\end{equation*} 
and the central limit theorem property:
for all $u$ H\"older continuous,
the random variable 
\begin{equation*}
\frac{1}{\sqrt{n}}\sum\limits_{j=0}^{n-1}(u\circ f^j - \int_M u d\mu) 
\end{equation*}
converges to a Gaussian law. More precisely, let $u$ be a $\alpha$-H\"older 
continuous function and  
\begin{equation*}
\sigma^2:=\int v^2 d\mu + 2\sum\limits_{j=1}^{\infty} v\cdot (v\circ f^j) d\mu,
\quad \text{ where } \quad v=u-\int u d\mu.
\end{equation*}
Then $\sigma<\infty$ and $\sigma=0$ iff $u=\varphi\circ f - \varphi$ for some
$\varphi\in L^1(\mu)$. If $\sigma>0$, then, for every interval $A\subset\real$, 
\begin{equation*}
\mu\left(x\in M: \frac{1}{\sqrt{n}}\sum\limits_{j=0}^{n-1}
\left(u(f^j(x)-\int u d\mu)\right)\in A\right)\to \frac{1}{\sigma\sqrt{2\pi}}
\int_A e^{-\frac{t^2}{2\sigma^2}} dt,
\end{equation*}
as $n\to\infty$.
\end{maintheorem} 

Let us comment about our conditions on the constants $\delta_0,\beta$ and
$\ep_0(f)$. We start with three conditions which should hold in the context of
both theorems: 
  
Most of the arguments of this paper strongly relies
on the non-uniform expansion property. A natural way to find it 
is to count the number of 
itineraries with frequent visit to $\SB= R_1\cup\dots, R_q$. 
Before doing this, 
we recall the following lemma
(which is a direct corollary of Stirling's formula).

\begin{lemma}\label{l.2.7}Given $\gamma\in (0,1)$, let $I_{\gamma,n}$ be
$$\left\{(i_0,\dots,i_{n-1})\in\{1,\dots,q,q+1,\dots,p+q\}^n; \ 
\#\{0\leq j < n; \ i_j\leq q\}>\gamma n\right\}$$
and $c_{\gamma}=\limsup\limits_{n\to\infty}\frac{1}{n}\log\# I_{\gamma,n}$.
Then, $c_{\gamma}\to\log q$ when $\gamma\to 1$.
\end{lemma} 

\begin{proof}See the appendix I for a proof. 
\end{proof}

We fix some $\gamma$ close to $1$ such that 
$c_{\gamma}<\log k$ (which is possible since $q<k$) and take $\delta_0$ 
small depending on $\sigma_1$ and $\gamma$ satisfying
\begin{equation}\label{e.4}
(1+\delta_0)^{\gamma} \sigma_1^{-(1-\gamma)}\leq e^{-4c}<1,
\end{equation}
for some constant $c>0$.\footnote{This condition ensures that if the orbit of a
point does not frequently visit $\SB$, then we have non-uniform expansion along
it.} 

Denote $c_0(f)=c_{\gamma}$ and consider $\ep_0(f)$ such that\footnote{Roughly 
speaking,~(\ref{e.epsilon2}) forces the orbit of generic points of certain
measures (e.g., equilibrium states) to stay a large portion of time outside 
$\SB$. See lemma~\ref{p.3.3}.} 
\begin{equation}\label{e.epsilon2}
\ep_0 (f)<\log k - c_0 (f).
\end{equation}

Also we assume that $\ep_0 (f)$ verifies
\begin{equation}\label{e.epsilon3}
e^{\ep_0(f)}\cdot\left(\frac{k-q}{k}\sigma_1^{-\alpha} + 
\frac{q}{k}(1+\delta_0)^{\alpha}\right)<1
\end{equation}
Clearly it is satisfied if $\ep_0(f)$ is small depending
only on $\sigma_1,k,q,\alpha$, since $k>q$.\footnote{The role 
of~(\ref{e.epsilon3}) is to guarantee that the equilibrium 
state obtained by our construction has exponential decay of correlations and 
the central limit theorem property.} 

Now we turn to the specific conditions of the main theorems:

\begin{itemize}
\item\textit{Precise conditions on the constants in theorem~\ref{t.A}.} 
Here, the hypothesis (H1), (H2), (H3) and the 
restrictions~(\ref{e.4}),~(\ref{e.epsilon2}),~(\ref{e.epsilon3}) are
\emph{almost} sufficient and, in fact, it is necessary only to complement the
restriction~(\ref{e.epsilon2}) with
\begin{equation}\label{epsilon2'}
\ep_0(f)<\log k - \max\limits_{x\in M}\log\|\Lambda^{d-1} Df(x)\|.
\end{equation} 

\item\textit{Precise conditions on the constants in theorem~\ref{t.B}.} Besides 
the hypothesis (H1), (H2), (H4) and the 
restrictions~(\ref{e.4}),~(\ref{e.epsilon2}),~(\ref{e.epsilon3})
above, we need two more
conditions. First of all, in the appendix
of~\cite{ABV}, Alves, Bonatti and Viana proved that there exists some
$\gamma_0>0$ depending only on $\sigma_2,p,q$ such that the orbit of Lebesgue
almost every point $x\in M$ spends a fraction $\gamma_0$ of the time in 
$R_{q+1}\cup\dots\cup R_{p+q}$:
$$\frac{1}{n}\#\{0\leq j
< n; \ f^j(x)\in R_{q+1}\cup\dots\cup R_{p+q}\}\leq\gamma_0 n,$$
for every large $n$. 

We take $\gamma_0<\gamma<1$ close to $1$ such that 
\begin{equation}\label{e.5}
\gamma m_2 + (1-\gamma) M_2 < \gamma_0 m_1 + (1-\gamma_0) M_1 - 
d\log (1+\de_0),
\end{equation}
where $d$ is the dimension of $M$. Note that by the hypothesis (H4), 
$m_2 < M_1$ and $m_2-m_1 < \beta$. Hence, this
condition can be verified if $\delta_0, \beta$ are sufficiently small
(depending only on $\gamma_0$, i.e., $\sigma_2,p,q$, and $M_1-m_2$). In 
particular, our condition ``$1-\gamma$ is small'' above depends 
only on $\sigma_2,p,q$ and $M_1-m_2$.

So, we get that $\beta$ is small depending only on $\sigma_2,p,q$, $M_1-m_2$ 
and $\delta_0$ is small depending on $\sigma_1,\sigma_2,p,q$, $M_1-m_2$.

Now we turn our attention to what does ``$\ep_0(f)$ is small'' in hypothesis
(H2) means. We take any 
$0<\rho$ such that 
$$\rho\cdot\left( \gamma_0 m_1 + (1-\gamma_0) M_1 - 
d\log (1+\de_0)\right)<\gamma m_2 + (1-\gamma) M_2.$$
Note that $\rho<1$ by~(\ref{e.5}). The last condition on $\ep_0(f)$ is 
\begin{equation}\label{e.epsilon1}
\ep_0(f)< (1-\rho)h_{top} (f).
\end{equation}

Finally, in the case of theorem~\ref{t.B}, we impose all the conditions above 
on the parameters $\beta,\delta_0$ and $\ep_0(f)$. 

For future reference, we recall that Oliveira proved 
\begin{theorem}[Oliveira~\cite{O}]\label{t.K}In the context of the
theorem~\ref{t.B}, i.e., under our hypothesis (H1), (H2), (H4), 
and if the previous conditions on $\beta,\delta_0$ and $\ep_0(f)$ holds, 
there are (ergodic) 
equilibrium states of $f$ for the potential $\phi$. Moreover, any 
ergodic equilibrium state $\eta$ of $f$ with respect to $\phi$ satisfies 
$\eta (V)\leq\gamma$. Furthermore, the Lyapounov exponents of $\eta$ are all
positive and $\eta$ has infinitely many hyperbolic times (see 
definition~\ref{d.hyp-times} below).
\end{theorem}
\end{itemize}

After this preparation, we are ready to study the Ruelle-Perron-Frobenius 
transfer operator and its spectral properties.

\section{The Ruelle-Perron-Frobenius operator}

The Ruelle-Perron-Frobenius transfer operator $\SL_\phi: C(M)
\rightarrow C(M)$ of $f:M\to M$ is defined on the space $C(M)$ of
continuous functions $g:M\to\real$ by
$$
\SL_\phi g(x) = \sum_{f(y)=x} e^{\phi(y)}g(y).
$$
Denote by $\lambda$ the spectral radius of $\SL_\phi$. Let $\nu$ be an 
eigenmeasure associated to the eigenvalue $\lambda$ of the dual operator 
$\SL_\phi^\star$ (i.e., $\SL_\phi^\star (\nu )= \lambda\cdot\nu$). In 
particular, $\lambda=\int\SL_\phi 1 d\nu$ and, hence, 
$k\leq \lambda\leq k\cdot e^{\max\phi}$. 

Note that the eigenmeasure $\nu$ is not necessarily $f$-invariant. However, we 
can capture its dynamical features (with respect to $f$) through the 
concept of Jacobian of a measure $\mu$ with respect to 
$f$: it is the (essentially)
unique function $J_{\mu}f$ such that
$$\mu (f(A)) = \int_A J_{\mu}f ,$$
for any measurable set $A$ where $f|_A$ is injective. Equivalently, $J_{\mu} f =
\frac{d ( f_{*}\mu )}{d\mu}$. In general, Jacobians do not exist but, if $f$ is
countable to one then $J_{\mu}f$ does exists for any measure $\mu$. 

As a first job, let us compute the Jacobian of the eigenmeasure $\nu$: 

\begin{lemma}\label{l.3.1}If $\SL_{\phi}^{\star}\nu = \lambda\nu$, then
$J_{\nu}f=\lambda e^{-\phi}$.
\end{lemma}

\begin{proof}Let $A$ be any measurable set such that $f|_A$ is injective. Take a
sequence $g_n$ of continuous functions satisfying $g_n\to\chi_A$ $\nu$-a.e. and
$\sup |g_n|\leq 2$ for all $n$. By definition,
$$\SL_{\phi} (e^{-\phi}g_n)=\sum\limits_{f(y)=x}g_n(y).$$
Since the right-hand of this equation is $\nu$-a.e. 
converging to $\chi_{f(A)}(x)$, applying the dominated convergence theorem and
the hypothesis $\SL_{\phi}^{\star}\nu = \lambda\nu$, we
get
$$\int\lambda e^{-\phi}g_n \ d\nu = \int e^{-\phi} g_n \
d(\SL_{\phi}^{\star}\nu) = \int\SL_{\phi} (e^{-\phi}g_n) \ d\nu\to \nu (f(A)).$$
However, the left-hand of the last equation converges to $\int_A \lambda
e^{-\phi} d\nu$. Thus,
$$\nu (f(A))=\int_A\lambda e^{-\phi} d\nu .$$
In other words, $J_{\nu}f=\lambda e^{-\phi}$.
\end{proof}

Next, we will show how the explicit formula for the Jacobian of $\nu$ can be
used to get some non-trivial information about the dynamics of $f$-orbits of
$\nu$ generic points.
  
\subsection{Non-uniform expansion and its consequences}

\textbf{Standing hypothesis of this subsection.} For the proof of the results of
the present subsection, we assume only (H1), (H2) and the 
conditions~(\ref{e.4}),~(\ref{e.epsilon2}).

\bigskip  

Note that~(\ref{e.epsilon2}) and lemma~\ref{l.3.1} lead us to 
\begin{equation*}
J_{\nu}f\geq \lambda e^{-\max\phi}\geq e^{\log k - \max\phi}> e^{c_0(f)}.
\end{equation*}
Hence, we can find some $\kappa > c_0(f)$ such that
\begin{equation}\label{e.epsilon2'}
J_{\nu}f > e^{\kappa} > q.
\end{equation}

From this lower bound on $J_{\nu} f$, we will derive a relevant combinatorial
result about the visit of $\nu$-generic $f$-orbits to the bad region $\SB$ of
eventual contraction: 
 
\begin{lemma}\label{p.3.3}Any probability measure $\nu$ such that 
$\SL_{\phi}^{\star}\nu = \lambda\nu$ satisfies 
$$\nu (G) = 1,$$
where $G:=\left\{x\in M : \ \limsup\limits_{n\to\infty} 
\frac{1}{n}\#\{0\leq j\leq
n;\ f^j(x)\in\SB\}\leq\gamma\right\}$.
\end{lemma}

\begin{proof}Consider any cylinder $R^n\in\SR^n$. Since $f$ is injective on each
element of $\SR$, it follows that $f^n$ is injective on each element of $\SR^n$.
Hence, by the definition of Jacobian and the inequality~(\ref{e.epsilon2'}),
$$1\geq\nu (f^n(R^n)) = \int_{R^n}J_{\nu}f^n d\nu =
\int_{R^n}\prod\limits_{j=0}^{n-1} (J_{\nu}f\circ f^j) d\nu \geq 
e^{\ka n}\nu (R^n).$$
In particular $\nu (R^n)\leq e^{-\ka n}$. 

Let $A(n)$ be the union of all
$R^n\in I_{\ga ,n}$, that is, all $n$-cylinders $R^n = [i_0,\dots , i_{n-1}]$
such that
$$\frac{1}{n}\#\{0\leq j<n;\ i_j\leq q\} > \ga.$$
Recall the terminology of the lemma~\ref{l.2.7}. The cardinality of $I_{\ga ,n}$
grows with $n$ at the exponential rate $c_{\ga}$. Note that we have $\ka >
c_{\ga}$. So, for some $C>0$, 
$$\nu (A(n)) = \sum\limits_{R^n\in I_{\ga ,n}}\nu (R^n) \leq 
C e^{c_{\ga}n} e^{-\ka n},$$
i.e., $\nu (A(n))\to 0$ exponentially fast. 
Then, by Borel-Cantelli lemma, for
every $x$ in a full $\nu$-measure set and for every $n$ large enough, we have
$x\notin A(n)$. This finishes the proof, since $x\notin A(n)$ means exactly that
$R^n(x)=[i_0,\dots ,i_{n-1}]$ satisfies $\frac{1}{n}\#\{0\leq j<n;\ i_j\leq q\}
\leq \ga$ and $i_j\leq q$ is equivalent to $f^j(x)\in\SB$.
\end{proof}

To convert this combinatorial result into useful information, we recall the
definition of hyperbolic times: 
 
\begin{definition}\label{d.hyp-times}We say that $n\in\natural$ is an 
$c$-\emph{hyperbolic time} 
for $x\in M$ if
\begin{equation*}
\prod\limits_{k=0}^{j-1}\|Df(f^{n-k}(x))^{-1}\|\leq e^{-2cj} \quad \text{ for
every } 1\leq j\leq n,
\end{equation*}
and $n\in\natural$ is an hyperbolic time for a cylinder $R^n\in\SR^n$ if $n$ is
a hyperbolic time for all $x\in R^n$. We denote by $\SR_h^n$ the set of cylinder
$R^n\in\SR^n$ with $n$ as an hyperbolic time.
\end{definition}

Hyperbolic times are interesting objects because of the following key 
properties:

\begin{lemma}\label{l.3.5}There exists $\delta >0$ depending only on $f$ and
$c$, such that given any hyperbolic time $n\geq 1$ for a point $x\in M$, and
given any $1\leq j\leq n$, there is an inverse branch $f_{x,n}^{-j}$ of $f^j$ 
defined on the whole ball of radius $\delta$ around $f^n(x)$, which 
sends $f^n(x)$ to $f^{n-j}(x)$ and satisfies 
$$ d(f_{x,n}^{-j}(z),f_{x,n}^{-j}(w))\leq e^{-jc}d(z,w)$$
for every $z,w$ in the ball $B(f^n(x),\delta)$.
\end{lemma}

\begin{proof}See~\cite{ABV}.
\end{proof}

\begin{corollary}\label{c.3.6}There exists $C>0$ depending on $\SR$ and
$\delta >0$ such that if $R^n\in\SR_h^n$, then 
$$d(f^{n-j} (x),f^{n-j} (y))\leq C e^{-jc} d(f^n (x),f^n (y)),$$
for every $x,y\in R^n$.
\end{corollary}

\begin{proof}It suffices to cover the atoms $R_1,\dots ,R_{p+q}$ with a finite
number of balls of radius $\delta$ and the claim follows from a direct
application of the lemma~\ref{l.3.5} and the fact that the inner diameter of
each $R_i$ is finite. 
\end{proof}

\begin{corollary}\label{c.3.7}Given any $\al$-H\"older continuous function 
$\phi : M\to M$, there exists a constant $A$ such that for all 
$R^n\in\SR_h^n$ and
$x,y\in R^n$ holds
$$|S_n\phi (x) - S_n\phi (y)|\leq A \ d(f^n(x),f^n(y))^{\al}.$$
\end{corollary}

\begin{proof}Since $\phi$ is $\al$-H\"older continuous,
$$|S_n\phi (x) - S_n\phi (y)|\leq 
\sum\limits_{j=0}^{n-1}|\phi (f^j(x)) - \phi (f^j(y))| \leq 
\|\phi\|_{\al}\sum\limits_{j=0}^{n-1}d(f^j(x),f^j(y))^{\al}.$$
By the corollary~\ref{c.3.6} we have 
$d(f^{n-j}(x),f^{n-j}(y))\leq Ce^{-jc}d(f^n(x),f^n(y))$ for every $j=0,\dots
,n-1$. These two facts together implies
$$|S_n\phi (x) - S_n\phi (y)|\leq 
\|\phi\|_{\al}C^{\al}\sum\limits_{j=0}^{n-1}e^{-jc\al}d(f^n(x),f^n(y))^{\al} 
\leq A \ d(f^n(x),f^n(y))^{\al}.$$
\end{proof}

\begin{corollary}\label{c.3.8}There exists $K$ such that if $R^n\in\SR_h^n$ and
$x,y\in R^n$ then 
$$K^{-1}\leq\frac{J_{\nu}f^n(x)}{J_{\nu}f^n(y)}\leq K.$$
\end{corollary}

\begin{proof}Notice that $J_{\nu}f^n(x)=\lambda^n e^{-S_n\phi(x)}$. Therefore,
the corollary~\ref{c.3.7} implies
$$\frac{J_{\nu}f^n(x)}{J_{\nu}f^n(y)} = e^{S_n\phi(y) - S_n\phi(x)} \leq 
e^{A \ d(f^n(x),f^n(y))^{\al}},$$
which gives the claimed result for $K = e^{AD^{\al}}$, where $D$ is the 
supremum of the diameters of the rectangles of $\SR$. 
\end{proof}

One of the main features of the measures considered here is contained in the
definition in the sequel: 

\begin{definition}We say that a measure $\eta$ is \emph{expanding} if
$\eta(H)=1$, where 
$$H:=\left\{x\in M : \ \limsup\limits_{n\to\infty} 
\frac{1}{n}\#\{0\leq j\leq
n;\ R^j(x)\in\SR_h^j\}>0\right\}.$$
\end{definition}

To study the expansiveness properties of measures we need the following
classical lemma:  

\begin{lemma}[Pliss~\cite{P}]\label{l.Pliss}Given $A\geq c_2 > c_1 > 0$, define
$\theta=\frac{c_2-c_1}{A-c_1}$. If $b_1,\dots,b_n$ are real numbers such that
$b_i\leq A$ and 
$$\sum\limits_{i=1}^{n}b_i\geq c_2 n,$$
then there are $l>\theta n$ integers, say $1\leq n_1<\dots<n_l\leq n$, such that
for $0\leq k\leq n_i$:
$$\sum\limits_{j=k+1}^{n_i}b_j\geq c_1 (n_i-k).$$
\end{lemma}

\begin{proof}Define $S(m)=\sum\limits_{j=1}^m (a_j - c_1)$, for $1\leq m\leq n$
and $S(0)=0$. Then, let $1< n_1 <\dots < n_l\leq n$ be the maximal sequence with
$S(n_i)\geq S(m)$ for all $0\leq m\leq n_i$, $i=1,\dots,l$. Since $S(n)>S(0)$,
$l>0$. On the other hand, the definitions means that 
\begin{equation*}
\sum\limits_{j=m+1}^{n_i}a_j\geq c_1 (n_i-m) \quad \text{ for } 0\leq m < n_i
\quad \text{ and } i=1,\dots,l.
\end{equation*}
Therefore, it remains only to show that $l>\theta n$. Observe that the defintion
of $S(m)$ implies
\begin{equation*}
S(n_i - 1) < S(n_{i-1}) \rightarrow S(n_i) < S(n_{i-1})+ A-c_1,
\end{equation*}
for $1<i\leq l$ and, moreover, $S(n_1)\leq A-c_1$ and $S(n_l)\geq S(n)\geq 
n(c_2 - c_1)$. Thus,
\begin{equation*}
n(c_2 - c_1)\leq S(n_l)=\sum\limits_{i=2}^l(S(n_i)-S(n_{i-1}))+S(n_1) <
l(A-c_1).
\end{equation*}
This completes the proof.
\end{proof}

At this point, the combinatorial result of lemma~\ref{p.3.3} permit us to show 
that any eigenmeasure of the dual of the transfer operator is expanding: 

\begin{lemma}\label{p.3.11}Let $\nu$ be a probability measure such that 
$\SL_{\phi}^{\star}\nu = \lambda\nu$. Then, there exists $\theta>0$ depending
only on $f$ and $c$ such that
$$\nu\left(\left\{x\in M : \ \limsup\limits_{n\to\infty} 
\frac{1}{n}\#\{0\leq j\leq
n;\ R^j(x)\in\SR_h^j\}\geq\theta\right\}\right) = 1.$$
In particular, $\nu$ is expanding.
\end{lemma}

\begin{proof}From the proof of the lemma~\ref{p.3.3}, we know that, if we define
$A(n)$ as the union of the $n$-cylinders $R^n = [i_0,\dots, i_{n-1}]$ such that 
$$\frac{1}{n}\#\{0\leq j < n; \ i_j\leq q\} >\gamma,$$
then the sequence $\nu (A(n))$ decreases exponentially with $n$. Take any
$x\notin A(n)$. Then, if $(i_0,\dots,i_{n-1})$ is the itinerary of $x$, 
$$\frac{1}{n}\sum\limits_{j=0}^{n-1}\log\sup\limits_{R_{i_j}}\|Df^{-1}\|\leq
\log\left( (1+\delta_0)^{\gamma}\sigma_1^{-(1-\gamma)}\right) \leq -4c <0.$$
Taking $A=\sup\limits_{x\in M}-\log\|Df(x)^{-1}\|$, $c_1=2c$, $c_2=3c$ and
$b_j=-\log\sup\limits_{R_{i_j}}\|Df^{-1}\|$ in lemma~\ref{l.Pliss}, we obtain
that there exists $n_1<\dots<n_l\leq n-1$ with $R^{n_i}\in\SR_h^{n_i}$ and
$l>\theta n$, where $\theta=\frac{c}{A-2c}$. These two facts finishes the proof.
\end{proof}

An interesting consequence of a measure $\eta$ being expanding is: 

\begin{lemma}\label{l.6.2}Let $\eta$ be an expanding measure. Then, for any 
measurable set 
$E$ and any $\ep>0$, there exists hyperbolic cylinders $C_1,\dots ,C_k$ such
that
$$\eta (E\Delta\bigcup\limits_{i=1}^k C_i) <\ep .$$ 
\end{lemma}

\begin{proof}Take $K_1\subset E$ and $K_2\subset E^c$ compact sets such that
$\eta (E\Delta K_1)<\ep / 3$, $\eta (E^c\Delta K_2)<\ep / 3$ and define
$r:=\text{dist}(K_1 , K_2 )$. Note that if $n$ is a hyperbolic time for some
cylinder $C$, then $\text{diam} (C)\leq K e^{-cn}$ (see corollary~\ref{c.3.6}). 
Hence, if $C\in\SR_h^n$ and $n>n_0$
(with $n_0$ large enough), we have $\text{diam} (C)\leq K e^{-cn}\leq r$. Since
$\eta$ is expanding (i.e., $\eta (G)=1$), we can choose $C_i\in\SR_h^{n_i}$
intersecting $K_1$ with
$n_i>n_0$ and 
$$\eta (K_1 \De \bigcup\limits_{i=1}^k C_i) < \ep / 3.$$
On the other hand, $C_i\cap K_2 =\emptyset$ because $C_i\cap K_1\neq\emptyset$,
$\text{diam} (C)\leq r$ and $r=\text{dist}(K_1 , K_2 )$.

Therefore,
$\eta (E\De\bigcup\limits_{i=1}^k C_i)\leq \eta (E-K_1) + \ep / 3 + \eta
(E^c-K_2)\leq\ep $.
\end{proof}

Another useful consequence of expansiveness is the \emph{weak Gibbs property}
for the eigenmeasures of $\SL_{\phi}^{\star}$:

\begin{lemma}\label{l.3.12}$\exists \ K>0$ such that 
if $\SL_{\phi}^{\star}\nu = \lambda\nu$, then for every $n$ and $x\in
R^n$ with $R^n\in\mathcal{R}^n_h$,

$$K^{-1}\leq\frac{\nu(R^n)}{\exp(S_n\phi(x) - Pn)}\leq K, $$
where $P=\log\lambda$. 
\end{lemma}

\begin{proof}By the lemma~\ref{l.3.1}, $J_{\eta}f^n = e^{nP-S_n\phi}$. Hence, if
$R^n = [i_0,\dots ,i_{n-1}]$, 
$$\eta (f(R_{i_{n-1}})) = \eta (f^n(R^n)) = \int_{R^n} J_{\eta}f^n d\eta =
\int_{R^n} e^{nP-S_n\phi} d\eta.$$
By the corollary~\ref{c.3.8}, there exists $K_1$ (not depending on $n$) 
such that for every $x,y\in R^n$
$$K_1^{-1} J_{\eta}f^n (y)\leq J_{\eta}f^n (x)\leq K_1 J_{\eta}f^n (y).$$
It follows that
$$K_1^{-1} \eta (f(R_{i_{n-1}}))\leq \frac{\eta(R^n)}{e^{S_n\phi - nP}}\leq K_1
\eta (f(R_{i_{n-1}})).$$
So, in order to conclude the proof, it suffices to verify that $\eta (f(R_i))>0$
for every $i$. Indeed, from the previous estimate, we have
$$K^{-1}\leq \frac{\eta(R^n)}{e^{S_n\phi - nP}}\leq K,$$
with $K=K_1\cdot\max\{\sup\limits_i \eta (f(R_i)), 
\sup\limits_i\frac{1}{\eta (f(R_i))}\}$.

By the third condition in the hypothesis (H1), we know that there
exists some $N$ satisfying $f^N(R_i)=M$ for every $i$. Since $\eta$ has a
Jacobian and $M$ has total probability, a decomposition of $f(R_i)$ into 
finitely many pieces where $f^{N-1}$ is injective shows that $\eta (f(R_i))>0$,
as desired.  
\end{proof}

An immediate corollary of the weak Gibbs property is the equivalence of any two
eigenmeasures of $\SL^{\star}$: 

\begin{corollary}\label{c.unique}Any two eigenmeasures $\eta_1$ and $\eta_2$ of
$\SL_{\phi}^{\star}$ with same eigenvalue $\lambda$ are equivalent.
\end{corollary}

\begin{proof}Using lemma~\ref{l.3.12}, we have for each $R\in\SR_h^n$,
$$ K^{-1} e^{S_n\phi (x) - Pn}\leq\nu_i(R)\leq K e^{S_n\phi (x) - Pn},$$
for $i=1,2$. This implies that $K^{-2}\nu_2(R)\leq\nu_1(R)\leq K^2\nu_2(R)$. The
proof is completed by applying lemma~\ref{l.6.2} and lemma~\ref{p.3.11}.
\end{proof}

The discussion of this subsection closes the study of the dual 
$\SL_{\phi}^{\star}$ of the transfer operator and its eigenmeasures (when the
local diffeomorphisms $f$ fits our assumptions). Now, we are going to prove 
some good spectral properties of the trasnfer operator $\SL_{\phi}$ itself. 

\subsection{Invariant cones for the transfer operator}\label{s.s.cones}

\textbf{Standing hypothesis of this subsection.} For the proof of the results of
the present subsection, we assume only (H1), (H2) and the 
condition~(\ref{e.epsilon3}).

\bigskip  

We define the norm $|||g|||_{\alpha} :=\sup\limits_{R_i}||g|_{R_i}||_{\alpha}$, where
$$||g||_{\alpha}:=\sup\limits_{x\neq y}\frac{|g(x)-g(y)|}{d(x,y)^{\alpha}} .$$

In order to prove the main results, we obtain the following Lasota-Yorke type 
inequality:

\begin{lemma}\label{l.LY}
$$|||\widetilde{\SL}_\phi(g)|||_{\alpha} \leq 
\Theta\cdot |||g|||_{\alpha}+C\cdot \|g\|_{\infty}, $$ where 
$\widetilde{\SL}_\phi :=\frac{1}{\lambda}\cdot\SL_\phi$ is the normalized 
Ruelle-Perron-Frobenius operator and $C>0$, $0<\Theta <1$ are real constants.
\end{lemma}

\begin{proof} For $x, x'$ in the same atom of the partition $\SR=\{R_1,\dots,R_d\}$, we have 
$$\frac{|\widetilde{\SL}_\phi g(x)-\widetilde{\SL}_\phi g(y)|}{d(x,x')^{\alpha}}\leq\frac{1}{\lambda} 
\sum\limits_{i}\frac{|e^{\phi(y_i)}g(y_i)-e^{\phi(y_i')}g(y_i')|}{d(x,x')^{\alpha}}\leq
\Theta\cdot |||g|||_{\alpha}+C\cdot \|g\|_{\infty}, $$
where $$\Theta=\frac{p(x)}{\lambda}e^{\max\phi}\sigma_1^{-\alpha} + 
\frac{q(x)}{\lambda}e^{\max\phi}(1+\delta_0)^{\alpha}$$ 
and $$C=\frac{|||e^{\max\phi}|||_{\alpha}}{\lambda}\big(p(x)\sigma^{-\alpha}+
q(x)(1+\delta_0)^{\alpha}\big).$$
Here $p(x)$ is the number of pre-images of $x$ belonging to
$\bigcup_{i=q+1}^{p+q}R_i$ and $q(x)$ is the number of pre-images of $x$ 
belonging to $\bigcup_{i=1}^{q}R_i$. This completes the proof because the
condition~(\ref{e.epsilon3}) implies $\Theta<1$. 
\end{proof}

The Lasota-Yorke type inequality above ensures the existence of a strictly 
$\SL_{\phi}$-invariant family of cones of positive H\"older continuous 
functions\footnote{C.M. would like to prof. Carlangelo Liverani for 
interesting discussions who lead the authors to the final 
expression for the cones $\Lambda_L$ used here.} In fact, if we put 
$\Lambda_L:=\{g\in C(M): g>0 \text{ and } |||g|||_{\alpha}\leq L \cdot 
\inf g\}$, it is not hard to prove that

\begin{corollary}\label{c.strict} 
There exists some $\sigma<1$ such that 
$\widetilde{\SL}_\phi(\Lambda_L)\subset\Lambda_{\sigma L}$, for
any $L$ sufficiently large (depending only on $\Theta$ and $C$).
\end{corollary}

\begin{proof} Let $g\in\Lambda_L$. By the lemma~\ref{l.LY} and the definition of
$\Lambda_L$,
$$|||\widetilde{\SL}_\phi g|||_{\alpha}\leq \Theta\cdot L\cdot \inf g +C\cdot
\|g\|_{\infty}.$$
Now, using that $\|g\|_{\infty}\leq \inf g +
d\vep^{\alpha} |||g|||_{\alpha}$, for any function $g$, we obtain
$$|||\widetilde{\SL}_\phi g|||_{\alpha}\leq (\Theta L + C_1) \inf g ,$$
where $C_1=C+C\cdot d\cdot\vep^{\alpha}$. However, we see that if $L$ is
sufficiently large (e.g., $L>\frac{C_1}{\Theta_0-\theta}$), then
$\Theta L + C_1 < \Theta_0 L$ (for any fixed $\Theta<\Theta_0<1$). 
In particular,
$$|||\widetilde{\SL}_\phi g|||_{\alpha}\leq \Theta_0\cdot L \cdot 
\inf g.$$
Finally, it is not difficult to see that $\inf g\leq
e^{\max\phi}\inf\widetilde{\SL}_\phi g$. So,
$$|||\widetilde{\SL}_\phi g|||_{\alpha}\leq \Theta_0\cdot e^{\max\phi}\cdot L 
\cdot\inf\widetilde{\SL}_\phi g.$$
Taking $\sigma=\Theta_0 \cdot e^{\max\phi}$ and using the
condition~(\ref{e.epsilon3}), we conclude the proof.
\end{proof}

Now, we compute the projective metric of $\Lambda_L$ 
(see the appendix II for precise definitions): 

\begin{lemma}\label{l.cone-metric}
The $\Lambda_L$-cone metric $\Psi_L$ is given by 
$\Psi_L(h,g)=\log\frac{B_L(h,g)}{A_L(h,g)}$, where
$$A_L(h,g)=\inf\limits_{x\neq y,\ x,y\in R_i\ z\in M } 
\frac{L|x-y|^{\al}g(z)-\left(g(x)-g(y)\right)}
{L|x-y|^{\al}h(z)-\left(h(x)-h(y)\right)} ,$$
and
$$B_L(h,g)=\sup\limits_{x\neq y,\ x,y\in R_i\ z\in M }
\frac{L|x-y|^{\al}g(z)-\left(g(x)-g(y)\right)}
{L|x-y|^{\al}h(z)-\left(h(x)-h(y)\right)} .$$
\end{lemma}
\begin{proof}By definition, $A h\preceq g$ iff 
$$g(x)-A h(x)\geq 0 \quad \forall x\in M,$$
and 
$$|||g-A h|||_{\al}\leq L\inf (g-A h).$$
I.e., 
$$A\leq\min\left\{ \inf\limits_{x} \frac{g(x)}{h(x)}, \ \ 
\inf\limits_{x\neq y,\ x,y\in R_i, \ z\in M}
\frac{L|x-y|^{\al}g(z)-\left(g(x)-g(y)\right)}
{L|x-y|^{\al}h(z)-\left(h(x)-h(y)\right)}\right\}.$$
But, if we choose $x_0\in M$ such that $\inf\limits_x \frac{g(x)}{h(x)} =
\frac{g(x_0)}{h(x_0)}$, then it is not hard to verify 
$$\frac{L|x-y|^{\al}g(x_0)-\left(g(x)-g(x_0)\right)}
{L|x-y|^{\al}h(x_0)-\left(h(x)-h(x_0)\right)}\leq
\frac{g(x_0)}{f(x_0)}.$$
Hence,
$$A_L(h,g)= \inf\limits_{x\neq y,\ x,y\in R_i, \ z\in M}
\frac{L|x-y|^{\al}g(z)-\left(g(x)-g(y)\right)}
{L|x-y|^{\al}h(z)-\left(h(x)-h(y)\right)}.$$
In a similar way,
$$B_L(h,g)=\sup\limits_{x\neq y,\ x,y\in R_i\ z\in M }
\frac{L|x-y|^{\al}g(z)-\left(g(x)-g(y)\right)}
{L|x-y|^{\al}h(z)-\left(h(x)-h(y)\right)} .$$
\end{proof}

The explicit formula of the projective metric of $\Lambda_L$ allow us to
prove that:

\begin{proposition}\label{p.diameter}$\Lambda_{\sigma L}$ has finite diameter 
in the $\Lambda_L$-cone metric $\Psi_L$:
\begin{equation*}
\text{diam}_{\Psi_L} (\Lambda_{\sigma L})\leq 
\Delta:=\Delta(\sigma,L):= 2\log\left( \frac{1+\sigma}{1-\sigma} \right) + 
2\log (1+\sigma L\cdot d\cdot\vep^{\al})
\end{equation*}
\end{proposition}

\begin{proof}
From the definition of $\Lambda_L$, we have for any $g\in\Lambda_{\sigma L}$,
$$|||g|||_{\al}\leq \sigma L\inf g.$$
In particular, we get 
$$\sup g\leq \inf g + \sigma L\cdot d\cdot\vep^{\al}\inf g.$$
However, the lemma~\ref{l.cone-metric} says that
$$
\Psi_L(h,g)=\log\sup\limits_{x,y\in R_i,\ u,v\in R_j,\ z,w\in M}
\frac{L\cdot g(z)-\frac{\left(g(x)-g(y)\right)}{|x-y|^{\al}}}
{L\cdot h(z)-\frac{\left(h(x)-h(y)\right)}{|x-y|^{\al}}}
\cdot\frac{L\cdot h(w)-\frac{\left(h(u)-h(v)\right)}{|u-v|^{\al}}}
{L\cdot g(w)-\frac{\left(g(u)-g(v)\right)}{|u-v|^{\al}}}.
$$
Putting these facts together, we
obtain
$$\Psi_L(h,g)\leq 2\log\left( \frac{1+\sigma}{1-\sigma} \right) + 
2\log (1+\sigma L\cdot
d\cdot\vep^{\al}) ,$$
for any $h,g\in\Lambda_{\sigma L}$. 
\end{proof}

From these results, it is quite standard to construct a 
$f$-invariant measure which is absolutely continuous with respect to $\nu$ 
whose density is H\"older continuous: 
 
\begin{theorem}\label{t.Ruelle1} There exists an unique $f$-invariant 
probability $\mu$ absolutely continuous with respect to $\nu$ whose 
Radon-Nikodym derivative $h$ is positive and H\"older continuous. 
\end{theorem}

\begin{proof}Before starting the proof, we show a simple lemma relating some 
Banach norms and the partial orderings induced by the cones $\Lambda_L$:
\begin{lemma}\label{l.Ruelle1}For all $\varphi,\psi\in\Lambda_L$, 
\begin{equation*}
-\varphi\preceq\psi\preceq\varphi \quad \text{ implies } \quad 
\begin{cases}
\int |\psi| d\nu\leq \int|\varphi| d\nu,\\
\sup\limits_M |\psi|\leq \sup\limits_M |\varphi| \text{ and }\\
|||\psi|||_{\alpha}\leq L\sup\limits_M |\varphi|.
\end{cases}
\end{equation*}
\end{lemma}

\begin{proof}[Proof of lemma~\ref{l.Ruelle1}]The relation
$-\varphi\preceq\psi\preceq\varphi$ means that $\varphi\pm\psi\in\Lambda_L$. In
particular, $\psi\leq\varphi$ and $-\psi\leq\varphi$. So, $\int |\psi| d\nu\leq
\int|\varphi| d\nu$ and $\sup\limits_M |\psi|\leq \sup\limits_M |\varphi|$.

To ensure that $|||\psi|||_{\alpha}\leq L\sup\limits_M |\varphi|$, we take
two points $x,y\in R_i$ in the same rectangle $R_i$ with $\psi(x)>\psi(y)$ and 
observe that
\begin{equation*}
\varphi (y)-\psi (y)\leq (1 + L d(x,y)^{\alpha})\cdot (\varphi (x)-\psi (x))
\end{equation*}
\begin{equation*}
\varphi (x)+\psi (x)\leq (1 + L d(x,y)^{\alpha})\cdot (\varphi (y)-\psi (y))
\end{equation*}
since $\varphi\pm\psi\in\Lambda_L$. Combining these inequalities, we obtain 
\begin{equation*}
\psi (x)-\psi (y) + \varphi (y) + \varphi (x) \leq (1+L d(x,y)^{\alpha})\cdot
 \left(\varphi (y) + \varphi (x) + \psi (y)-\psi (x)\right). 
\end{equation*}
Hence,
\begin{equation*}
\left(2+L d(x,y)^{\alpha}\right)\cdot (\psi (x)-\psi (y)) \leq L d(x,y)^{\alpha}
 (\varphi (x)+\varphi (y)).
\end{equation*}
Therefore,
\begin{equation*}
|||\psi|||_{\alpha}\leq L\sup\limits_M\varphi.
\end{equation*}
\end{proof}
Now, the corollary~\ref{c.strict} says: there exists $\sigma<1$ such that, for 
all $L>0$ sufficiently large, 
$\widetilde{\SL}_{\phi}(\Lambda_L)\subset \Lambda_{\sigma L}$. Applying the
theorem~\ref{t.Birkhoff} (see appendix II below) and
propostion~\ref{p.diameter}, we obtain, for arbitrary
$\varphi,\psi\in\Lambda_L$,
\begin{equation}\label{e.strict}
\begin{split}
\Psi_L(\widetilde{\SL}_{\phi}^n(\varphi), \widetilde{\SL}_{\phi}^{k+n}(\psi))
&\leq \left(\tanh(\Delta / 4)\right)^{n-1}
\Psi_L(\widetilde{\SL}_{\phi}(\varphi),\widetilde{\SL}_{\phi}^k(\psi)) \\
\left(\tanh(\Delta / 4)\right)^{n-1} \Delta,
\end{split}
\end{equation}
where $\Delta<\infty$ is the diameter of $\Lambda_{\sigma L}$ with respect to
the $\Lambda_L$-cone metric. Since 
$\int\widetilde{\SL}_{\phi}^n\varphi d\nu = \int\varphi d\nu$, a simple
consequence of the lemma~\ref{l.Ruelle1} (with $\|.\|_{(1)}:=\int |.| d\nu$ and 
$\|.\|_{(2)}:=
\sup\limits_M |.|+\frac{1}{L}|||.|||_{\alpha}$) and lemma~\ref{l.appendixII} (of 
appendix II below) is: $\widetilde{\SL}_{\phi}^n 1$
is a Cauchy sequence for the supremum norm. Its limit $h=\lim\limits_{n\to
\infty}\widetilde{\SL}_{\phi}^n 1$ is a fixed point of $\widetilde{\SL}_{\phi}$
and $h\in\Lambda_L$. This implies that $\mu = h \nu$ is an invariant
probability. Indeed, for a test function $u$,
\begin{equation*}
\int u\circ f d\mu = \lambda^{-1}\int\SL_{\phi}(h\cdot u\circ f) d\nu =
\lambda^{-1}\int\SL_{\phi}h\cdot u d\nu = \int h\cdot u d\nu = \int u d\mu.
\end{equation*} 
The uniqueness part of the statement of the theorem follows by taking a fixed
point $\varphi\in\Lambda_L$ with $\int\varphi d\nu = 1$ and using the
estimate~(\ref{e.strict}). 
\end{proof}

Finally, the Lasota-Yorke type inequality of lemma~\ref{l.LY} implies that the
$f$-invariant probability $\mu$ just obtained has fast decay of correlations: 

\begin{theorem}\label{t.decay}
$\mu$ has exponential decay of correlations for H\"older
observables: there exists some constants $0<\tau<1$ and $K>0$ such that,
for all $u\in L^1(\nu), v\in C^{\alpha}(M)$ and $n\geq 0$,
\begin{equation*}
\left|\int_M (u\circ f^n)v d\mu - \int_M u d\mu\int_M v d\mu\right|\leq 
K(u,v)\cdot\tau^n,
\end{equation*}
and $\mu$ satisfies the central limit theorem: for all $u$ H\"older continuous,
the random variable 
\begin{equation*}
\frac{1}{\sqrt{n}}\sum\limits_{j=0}^{n-1}(u\circ f^j - \int_M u d\mu) 
\end{equation*}
converges to a Gaussian law. More precisely, let $u$ be a $\alpha$-H\"older 
continuous function and  
\begin{equation*}
\sigma^2:=\int v^2 d\mu + 2\sum\limits_{j=1}^{\infty} v\cdot (v\circ f^j) d\mu,
\quad \text{ where } \quad v=u-\int u d\mu.
\end{equation*}
Then $\sigma<\infty$ and $\sigma=0$ iff $u=\varphi\circ f - \varphi$ for some
$\varphi\in L^1(\mu)$. If $\sigma>0$, then, for every interval $A\subset\real$, 
\begin{equation*}
\mu\left(x\in M: \frac{1}{\sqrt{n}}\sum\limits_{j=0}^{n-1}
\left(u(f^j(x)-\int u d\mu)\right)\in A\right)\to \frac{1}{\sigma\sqrt{2\pi}}
\int_A e^{-\frac{t^2}{2\sigma^2}} dt,
\end{equation*}
as $n\to\infty$.
\end{theorem}

\begin{proof} The plan of the proof is, firstly, to prove the exponential decay
of correlations and, then, to use this information to get the central limit 
theorem.

We note that the correlation function 
\begin{equation*}
C_{u,v} (n):=\int_M (u\circ f^n)v d\mu - \int_M u d\mu\int_M v d\mu.
\end{equation*}
can be rewritten as 
\begin{equation*}
C_{u,v} (n)=\int_M (u\circ f^n)v h d\nu - \int_M u d\mu\int_M v h d\nu.
\end{equation*}
Suppose first that $vh\in\Lambda_L$, for some sufficiently large $L$. It is no 
restriction to assume that $\int v h d\nu = 1$. Then, 
denoting $\|u\|_1:=\int |u| d\mu$,
\begin{equation*}
\begin{split}
\left|\int (u\circ f^n)v h d\nu - \int u d\mu\int_M v h d\nu\right| &= 
\left|\int u\left(\frac{\widetilde{\SL}_{\phi}^n (v h)}{h}-1\right)
d\mu\right|\\ 
&\leq \left\|\frac{\widetilde{\SL}_{\phi}^n(v h)}{h}-1\right\|_{0}\cdot\|u\|_1.
\end{split}
\end{equation*}
To estimate the term 
$\left\|\frac{\widetilde{\SL}_{\phi}^n(v h)}{h}-1\right\|_{0}$, we
note that $h>0$ is a (H\"older) continuous function on a compact manifold $M$; 
so, applying the estimate~(\ref{e.strict}) and the 
lemma~\ref{l.appendixII} of appendix II (with $\|.\|_{(1)}:=\int |.| d\nu$ and 
$\|.\|_{(2)}:= \sup\limits_M |.|+\frac{1}{L}|||.|||_{\alpha}$), we have 
\begin{equation*}
\|\frac{\widetilde{\SL}_{\phi}^n(v h)}{h}-1\|_{0}\leq (e^{D\tau^n}-1)\leq
K\tau^n,
\end{equation*}
for some constants $0<\tau<1$, $D=D(u,v)>0$ and $K=K(u,v)$.

In particular,
\begin{equation*}
\left|\int (u\circ f^n)v h d\nu - \int u d\mu\int_M v h d\nu\right| \leq 
K(u,v)\tau^n\|u\|_1,
\end{equation*}
if $vh\in\Lambda_L$. 

For a general $(A,\alpha)$-H\"older continuous function $v$, we write, for
$B>0$, $vh=\varphi$ and 
\begin{equation*}
\varphi = \varphi_B^+ - \varphi_B^- \quad \text{ where } \quad \varphi_B^+ =
\frac{1}{2} (|\varphi|\pm\varphi)+B.
\end{equation*}  

Clearly $\varphi_B^{\pm}$ are $(A,\alpha)$-H\"older continuous and
$\varphi_B^{\pm}\geq B$. Taking $B$ sufficiently large, we get
$\varphi_B^{\pm}\in\Lambda_L$, for any sufficiently large fixed $L$. So, we can 
apply the previous estimate for $\varphi_B^{\pm}$. Hence, by 
linearity, the same estimate holds for $\varphi$. Note that the constant
$K(u,v)$ has the form 
\begin{equation*}
K(u,v)\leq K_0\|u\|_1 (\|v\|_1+|||v|||_{\alpha}),
\end{equation*}
where $K_0$ is independent of $u,v$. In fact, the first part of the proof gives
\begin{equation*}
K(\varphi_B^{\pm},v)\leq 
K_1\|u\|_1\left(\|\varphi\|_1+B\right)\leq K_2\|u\|_1(\|v\|_1+|||v|||_{\alpha}),
\end{equation*}
and the claim follows since $K(u,v)\leq
K_3(K(u,\varphi_B^+)+K(u,\varphi_B^-)$.\footnote{This technical remark is useful
when proving the ergodicity of $\mu$, since it allows us to conclude uniform
bounds for the decay of the correlation function when the function $u$ varies
on a $L^1(\mu)$ bounded set. See the proof of lemma~\ref{l.exact} below.}

This concludes the proof of the exponential decay of correlations. Now we are
going to use this information to show the central limit theorem. But, before
performing the argument, we introduce some notation. Let $\SF$ be the
Borel sigma-algebra of $M$ and $\SF_n:=f^{-n}(\SF)$. Recall that a function
$\xi:M\to\real$ is $\SF_n$-measurable iff $\xi=\xi_n\circ f^n$ for some ($\SF$-,
i.e., Borel) measurable $\xi_n$. Note that $\SF_{i}\supset\SF_{i+1}$ for every
$i\geq 0$. Also, an $f$-invariant measure $\mu$ is \emph{exact} iff the
sigma-algebra 
\begin{equation*}
\SF_{\infty}=\bigcap\limits_{n\geq 0}\SF_n
\end{equation*} 
is $\mu$-trivial.\footnote{Equivalently, all $\SF_{\infty}$-measurable functions
are constant $\mu$-a.e.} Clearly any exact measure is ergodic, since
$\chi_A\in\SF_{\infty}$ if $A$ is $f$-invariant.

Coming back to the proof of the central limit theorem, a direct corollary of 
the exponential decay of correlations is the exactness of $\mu$: 

\begin{lemma}\label{l.exact} $\mu$ is exact.
\end{lemma} 

\begin{proof} Let $\psi\in L^1(\mu)$ be $\SF_{\infty}$-measurable. By
definition, there are $\psi_n$ Borel measurable functions with $\psi=\psi_n\circ
f^n$. Obviously, $\|\psi_n\|_1=\|\psi\|_1<\infty$. Therefore, our previous
discussion concernig the asymptotics of the correlation function allows us to
obtain, for any $\alpha$-H\"older continuous function $\varphi$, 
\begin{equation*}
\begin{split}
\left|\int(\psi-\int\psi d\mu)\varphi d\nu\right| &= 
\left|\int(\psi_n\circ f^n)\varphi-\int\psi d\mu\int\varphi d\nu\right| \\ 
&\leq K_0^{'}(\varphi)\|\psi_n\|_1\tau^n = K_0^{'}(\varphi)\|\psi\|_1\tau^n
\to 0,
\end{split}
\end{equation*}
as $n\to\infty$, for some $K_0^{'}(\varphi)>0$. So, 
$\int(\psi-\int\psi d\mu)\varphi d\nu = 0$ and, in particular, 
$\psi=\int\psi d\mu$ almost everywhere. In other words, $\mu$ is exact.
\end{proof}

Let $L^2(\SF_n)=\{\xi\in L^2(\mu):\xi \text{ is }\SF_n-\text{ measurable }\}$. 
Observe that $L^2(\SF_n)\supset L^2(\SF_{n+1})$ for each $n\geq 0$. Given
$\varphi\in L^2(\mu)$, we denote $\mathbb{E}(\varphi|\SF_n)$ the 
$L^2$-orthogonal projection of $\varphi$ to $L^2(\SF_n)$. 

Another easy consequence of the exponential decay of correlations is:

\begin{lemma}\label{l.projections}For every $\alpha$-H\"older continuous
function $u$ with $\int u d\mu = 0$ there is $R_0=R_0(u)$ such that
$\|\mathbb{E}(u|\SF_n)\|_2\leq R_0\tau^n$ for all $n\geq 0$.
\end{lemma} 

\begin{proof}
\begin{equation*}
\begin{split}
\|\mathbb{E}(u|\SF_n)\|_2 &= \sup\left\{\int\xi u d\mu: \xi\in L^2(\SF_n),
\|\xi\|_2 = 1\right\} \\
&= \sup\left\{\int(\psi\circ f^n) u h d\nu: \psi\in L^2(\mu),
\|\psi\|_2 = 1\right\} \\
&\leq K_0^{'}(uh)\tau^n,
\end{split}
\end{equation*}
since $\|\psi\|_1\leq\|\psi\|_2$ and $\int u d\mu = \int u h d\nu = 0$.
\end{proof}

Now the proof of the central limit theorem can be derived from the
lemmas~\ref{l.exact},~\ref{l.projections} and the following abstract lemma
(whose proof can be found in~\cite[p. 28--33]{V}):

\begin{lemma} Let $(M,\SF,\mu)$ be a probability space, $f:M\to M$ be a
measurable map such that $\mu$ is $f$-invariant and ergodic. Consider $u\in
L^2(\mu)$ such that $\int u d\mu = 0$ and denote by $\SF_n$ the non-increasing
sequence of sigma-algebras $\SF_n=f^{-n}(\SF), n\geq 0$. Assume that 
\begin{equation*}
\sum\limits_{n=0}^{\infty}\|\mathbb{E} (u|\SF_n)\|_2 <\infty.
\end{equation*}  
Then, the number $\sigma\geq 0$ defined by
\begin{equation*}
\sigma^2 = \int u^2 d\mu + 2\sum\limits_{n=1}^{\infty} u\cdot (u\circ f^n) d\mu
\end{equation*}
is finite and $\sigma = 0$ iff $u=\varphi\circ f - \varphi$ for some $\varphi\in
L^2 (\mu)$. On the other hand, if $\sigma>0$, then, for any interval
$A\subset\real$,
\begin{equation*}
\mu\left(x\in M: \frac{1}{\sqrt{n}}\sum\limits_{j=0}^{n-1}
\left(u(f^j(x)-\int u d\mu)\right)\in A\right)\to \frac{1}{\sigma\sqrt{2\pi}}
\int_A e^{-\frac{t^2}{2\sigma^2}} dt,
\end{equation*}
as $n\to\infty$.
\end{lemma}
This completes the proof of the theorem.
\end{proof}

At this stage, we can prove that the $f$-invariant measure $\mu$ constructed
above is an equilibrium state for $(f,\phi)$. 

\section{Existence of the equilibrium states}

We denote by $g:M\to (0,\infty)$ the function
\begin{equation*}
g(x):=\lambda^{-1} e^{\phi(x)}\frac{h(x)}{h(f(x))},
\end{equation*}
where $h>0$ is the eigenfunction of $\SL_{\phi}$ with eigenvalue $\lambda$. In
particular, we have, for every $x\in M$,
\begin{equation}\label{e.12}
\sum\limits_{f(y)=x} g(y) = 
\frac{\sum\limits_{f(y)=x}e^{\phi(y)}h(y)}{\lambda h(x)} = 
\frac{\SL_{\phi}h(x)}{\lambda h(x)} = 1.
\end{equation}

For later reference, we state the following elementary calculus lemma:

\begin{lemma}\label{l.6.4}
Let $p_i, x_i$, $i=1,\dots ,n$ be positive real numbers such that
$\sum\limits_{i=1}^n p_i =1$. Then, 
\begin{equation*}
\sum\limits_{i=1}^n p_i\log x_i \leq \log (\sum\limits_{i=1}^n p_i x_i),
\end{equation*}
and the equality holds if and only if the numbers $x_i$ are all equal.
\end{lemma}

We are ready to prove the main result of this section, namely the existence of
an equilibrium measure: 

\begin{proposition}\label{l.6.5}$\mu$ is an equilibrium state of $\phi$, 
$P=\log\lambda$ is the pressure
$P(f,\phi)$. Moreover, if $\eta$ is any equilibrium state of $\phi$,
then 
\begin{itemize}
\item $h_{\eta}(f)+\int\log g \ d\eta = 0$;
\item $J_{\eta}f(y)=1/g(y)$ for $\eta$-a.e. $y\in M$,
\end{itemize}
either in the case of theorem~\ref{t.A} or in the case of thereom~\ref{t.B}. 
\end{proposition}

\begin{proof}
[Proof of proposition~\ref{l.6.5} in the context of theorem~\ref{t.A}] 
Since $\mu$ is expanding, by the lemma~\ref{l.6.2} we know that
$\SR$ is a generating partition for $\mu$. Kolmogorov-Sinai's theorem implies 
$h_{\mu}(f)=h_{\mu}(f,\SR)$. The equality $h_{\mu}(f)+\int\phi \ d\mu =
\log\lambda = P$ follows from the lemma~\ref{l.3.12}.

Now let $\eta$ be an invariant measure such that $h_{\eta} (f)
+ \int\phi d\eta\geq\log\lambda$, and hence 
$$ h_{\eta} (f) +\int\log g \ d\eta = h_{\eta} (f) -\log\lambda + 
\int\left(\phi+\log h - \log h\circ f\right) d\eta\geq 0.$$
Applying Rokhlin's formula
(see appendix III, corollary~\ref{c.rokhlinA} below)\footnote{This is the unique
point where the extra hypothesis (H3) of theorem~\ref{t.A} is used. Indeed, if
we know that, \emph{a priori}, Rokhlin's formula is true for the measures $\eta$
with $h_{\eta}(f)+\int\phi d\eta\geq\lambda$, then (H3) could be removed.}
$$ h_{\eta} (f) = \int \log J_{\eta} f ,$$
into the previous estimate, we get 
$$0\leq\int\frac{g}{g_{\eta}} d\eta = \int\sum\limits_{f(y)=x} g_{\eta} (y)
\log\frac{g(y)}{g_{\eta}(y)} d\eta,$$
where $g_{\eta} := 1 / J_{\eta} f$. However, in view of~(\ref{e.12}) and 
lemma~\ref{l.6.4}, we get
$$ \sum\limits_{f(y)=x} g_{\eta} (y)\log\frac{g(y)}{g_{\eta}(y)}\leq
\log\left(\sum\limits_{f(y)=x}g_{\eta} (y)\frac{g(y)}{g_{\eta}(y)}\right) = 
\log (\sum\limits_{f(y)=x} g(y))=0 $$
at $\eta$-a.e. $x$. From this discussion, we conclude that 
$\sum\limits_{f(y)=x} g_{\eta} (y)\log\frac{g(y)}{g_{\eta}(y)}$ is a 
non-positive function with non-negative integral (with respect to $\eta$). 
Hence it vanishes at $\eta$-a.e. $x$. Thus, we have the equality 
$P(f,\phi) = h_{\eta}(f) + \int\phi d\eta = \log\lambda$. From the
lemma~\ref{l.6.4}, we also obtain that the values 
$\log\frac{g(y)}{g_{\eta}(y)}$ coincides for all $y\in f^{-1}(x)$. I.e., for
$\eta$-a.e. $x\in M$, there exists a number $c(x)$ such that
$$\frac{g(y)}{g_{\eta}(y)} = c(x) \quad \text{for every} \quad 
y\in f^{-1} (x).$$
Since $\eta$ is an invariant measure, 
$$\sum\limits_{f(y)=x}g_{\eta} (y) = 1,$$
for $\eta$-a.e. $x$. These two facts together with~(\ref{e.12}) gives
$$ c(x)=\frac{\sum\limits_{f(y)=x}g_{\eta} (y)}{\sum\limits_{f(y)=x}g (y)}=1 ,$$
at $\eta$-a.e. $x$. This implies $g(y)=g_{\eta}(y)$ for every $y$ on the
pre-image of a set of full $\eta$ measure. By invariance of $\eta$, this set
also has full $\eta$ measure. This completes the proof.
\end{proof}


\begin{proof}
[Proof of proposition~\ref{l.6.5} in the context of theorem~\ref{t.B}] 
The proof in this case is very similar to the previous one, except for the
class of invariant measures satisfying the Roklhin's formula. Again, since 
$\mu$ is expanding, by the 
lemma~\ref{l.6.2} we know that
$\SR$ is a generating partition for $\mu$. Kolmogorov-Sinai's theorem implies 
$h_{\mu}(f)=h_{\mu}(f,\SR)$. The equality $h_{\mu}(f)+\int\phi \ d\mu =
\log\lambda = P$ follows from the lemma~\ref{l.3.12}.

Now let $\eta$ be an ergodic equilibrium measure. In particular, $h_{\eta} (f)
+ \int\phi d\eta\geq\lambda$, and hence 
$$ h_{\eta} (f) +\int\log g \ d\eta = h_{\eta} (f) -\log\lambda + 
\int\left(\phi+\log h - \log h\circ f\right) d\eta\geq 0.$$
Proceeding as before, an application of Rokhlin's formula
(see appendix III, corollary~\ref{c.rokhlinB} below)\footnote{As before, this 
is also the unique
point where the extra hypothesis (H4) of theorem~\ref{t.B} is used. In fact, if
we know that, \emph{a priori}, Rokhlin's formula is true for ergodic equilibrium
measures $\eta$, then (H4) could be removed (as well as the
restrictions~(\ref{e.5}) and~(\ref{e.epsilon1}) ).}
$$ h_{\eta} (f) = \int \log J_{\eta} f $$
into the previous estimate still gives 
$$\frac{g(y)}{g_{\eta}(y)} = c(x) \quad \text{for every} \quad 
y\in f^{-1} (x).$$
Since $\eta$ is an invariant measure, 
$$\sum\limits_{f(y)=x}g_{\eta} (y) = 1,$$
for $\eta$-a.e. $x$. These two facts together with~(\ref{e.12}) gives
$$ c(x)=\frac{\sum\limits_{f(y)=x}g_{\eta} (y)}{\sum\limits_{f(y)=x}g (y)}=1 ,$$
at $\eta$-a.e. $x$. This implies $g(y)=g_{\eta}(y)$ for every $y$ on the
pre-image of a set of full $\eta$ measure. By invariance of $\eta$, this set
also has full $\eta$ measure. This completes the proof.
\end{proof}

To end this section, we show that equilibrium states are closely related to the 
eigenmeasures of $\SL^{\star}$:
 
\begin{lemma}\label{l.6.6}If $\eta$ is an equilibrium measure of $\phi$, then
$\SL_{\phi}^{\star}(h^{-1}\eta)=\lambda (h^{-1}\eta)$.
\end{lemma}

\begin{proof} Given any continuous function $\xi$,
$$\int\xi d(\SL_{\phi}^{\star} (h^{-1}\eta)) = \int h^{-1}(x)\SL_{\phi}\xi(x)
d\eta = \int\sum\limits_{f(y)=x} e^{\phi(y)}\cdot\frac{\xi(y)}{h(f(y))} 
d\eta (x).$$
From the definition of $g$ and the proposition~\ref{l.6.5}, 
$$ e^{\phi (y)} \frac{1}{h(f(y))}=\lambda\frac{g(y)}{h(y)} =
\lambda\frac{g_{\eta}(y)}{h(y)}.$$
Combining these two equations, we have 
$$\int\xi d(\SL_{\phi}^{\star} (h^{-1}\eta)) = \lambda
\int\frac{g_{\eta}(y)}{h(y)}\xi (y) d\eta (x) = \lambda\int\xi h^{-1} d\eta.$$
Because $\xi$ is arbitrary, it holds 
$\SL_{\phi}^{\star} (h^{-1}\eta)=\lambda (h^{-1}\eta)$.
\end{proof}

\section{Proof of the theorems~\ref{t.A} and~\ref{t.B}}

Once the previous results about the properties of the Ruelle-Perron-Frobenius
transfer operator are established, it is now an easy matter to put the
proposition~\ref{l.6.5}, lemma~\ref{l.6.6} and theorem~\ref{t.decay} together to
get the proof of our main theorems. 
   
The proposition~\ref{l.6.5} and theorem~\ref{t.decay} says that the 
measure $\mu$ constructed in the subsection~\ref{s.s.cones} is an equilibrium 
measure with exponential decay of correlations and the central limit theorem 
property (in both the contexts of theorem~\ref{t.A} and~\ref{t.B}). So, it 
remains only to prove that $\mu$ is the unique equilibrium state of $f$. 

By the lemma~\ref{l.6.6}, any ergodic equilibrium measure $\eta$ satisfies
$\SL^{\star}_{\phi}(h^{-1}\eta)=\lambda (h^{-1}\eta)$. Thus, if we denote
$\nu_{\eta}=h^{-1}\eta$, then $\nu_{\eta}$ is an eigenmeasure of
$\SL^{\star}_{\phi}$ with eigenvalue $\lambda$. However, the
corollary~\ref{c.unique} implies that $\mu$ and $\eta$ are equivalent measures,
that is, $\eta=\xi\mu$ for some $\mu$-integrable function $\xi$. Since $\eta$
and $\mu$ are invariant measures, $\eta = f^*\eta = (\xi\circ f)\mu$. It follows
that $\xi\circ f = \xi$ and, by the ergodicity of $\eta$, we get that $\xi$ is
constant. Using that $\eta$ and $\mu$ are probability measures, we obtain that
$\xi=1$, i.e., $\eta=\mu$. So, $\mu$ is the unique equilibrium
state of $f$. 

This completes the proof of theorems~\ref{t.A} and~\ref{t.B}.

\section{Final comments}

As a matter of fact, the machinery developed in the previous sections can be
applied to obtain further nice ergodic properties of the equilibrium states just
constructed. To illustrate the power of these tools, we take this last section
to do three remarks: 

\begin{itemize}
\item First, we show that the transfer operator can be approximated by finite
rank operator; thus, the transfer operator has the spectral gap property and
this produces a new proof of the exponential decay of correlations for the
equilibrium states considered here.
\item Second, we prove that the estimates of Lasota-Yorke type in
lemma~\ref{l.LY} are robust under small \emph{random perturbations};
consequently, using the work of Baladi and Young, we get that our equilibrium
measures are \emph{stochastically stable}. 
\item Third, we prove that the abundance (positive density) of hyperbolic times
for generic points (of the equilibrium states in theorems~\ref{t.A},~\ref{t.B})
and the uniqueness of the equilibrium measures $\mu_{\phi}$ are sufficient to 
conclude, along
the lines of the work of Young, that $\mu_{\phi}$ has the \emph{large deviation 
property}. 
\end{itemize}

Since these three remarks are simple modifications (based in the lemmas proved
here) of the standard arguments of the thermodynamical formalism of uniformly
expanding maps, we only sketch the proof of the claimed results by pointing out 
the relevant changes of the well-known arguments of the literature. 
    
\subsection{Spectral Gap}

The outline of this subsection is: we start with a new construction of the
equilibrium measures for the robust classes of maps (and potentials) considered
above; in the sequel, we show that these equilibrium measures are exact;
finally, we approximate the transfer operator by compact operators. The
conclusion is a new proof of the exponential mixing of these equilibrium states.

It is clear that the construction of the equilibrium states depends only on the
existence of an eigenfunction $h$ for the transfer operator $\SL_{\phi}$ with
eigenvalue $\lambda$. To accoplish this objective, we closely follow the
approach of Oliveira and Viana~\cite{OV2}: take the sequence of
functions
\begin{equation*}
g_n(x) = \sum\limits_{f^n(y)=x; \, R^n(y)\in\SR_h^n} e^{S_n\phi(y)}
\end{equation*} 
and the sequence of numbers 
\begin{equation*}
Z_n = \sum\limits_{R^n\in \SR^n_h} e^{S_n\phi(R^n)},
\end{equation*}
where $S_n\phi(R^n):=\max\{S_n\phi(y): y\in R^n\}$. Note that $g_n\leq Z_n$.

\begin{lemma} There exists a constant $K_1>0$ such that 
\begin{equation*}
\lambda^{-n}Z_n\leq K_1 \quad \textrm{ and } \quad K_1^{-1}\leq
\frac{1}{n}\sum\limits_{j=0}^{n-1} \lambda^{-j}Z_j\leq K_1,
\end{equation*}
for all $n\geq 0$.
\end{lemma}

\begin{proof}Since the reference measure $\nu$ is a non-lacunary weak Gibbs
measure (see lemma~\ref{l.3.12}), we know that 
\begin{equation*}
K^{-1}\nu(R^n)\leq\lambda^{-n}e^{S_n\phi(y)}\leq K\nu(R^n), 
\end{equation*}
for every $y\in R^n$, $R^n\in\SR_h^n$. So, if $H_n$ is the union of all
$R^n\in\SR_h^n$, it follows 
\begin{equation*}
K^{-1}\nu(H_n)\leq \lambda^{-n}Z_n\leq K\nu(B_n)
\end{equation*}
This proves the first part of the lemma $\lambda^{-n}Z_n\leq K$ since $\nu$ is a
probability measure. For the second part, note that we have 
\begin{equation*}
K^{-1}\frac{1}{n}\sum\limits_{j=0}^{n-1}\nu(B_j)\leq
\frac{1}{n}\sum\limits_{j=0}^{n-1}\lambda^{-j}Z_j.
\end{equation*}
Therefore, it suffices to obtain an uniform lower bound on the left hand side
quantity, which can be written as 
\begin{equation*}
\frac{1}{n}\sum\limits_{j=0}^{n-1}\nu(B_j) = \int\int\chi_{B_j}(x) d\nu(x)
dm_n(j),
\end{equation*}
where $m_n$ is the normalized counting measure of $\{1,\dots,n\}$. Because $\nu$
generic points have a positive density $\theta>0$ of hyperbolic times (see
lemma~\ref{p.3.11}), we have 
\begin{equation*}
\frac{1}{n}\sum\limits_{j=0}^{n-1}\chi_{B_j}(x)>\theta
\end{equation*}
for $\nu$ almost every $x\in M$ and all $n$ large enough (depending on $x$). 
Taking $n$ large so that the set $X$ of points $x\in M$ verifying this 
inequality has $\nu$ measure bigger than $1/2$, we can apply Fubini theorem to 
conclude
\begin{equation*}
\frac{1}{n}\sum\limits_{j=0}^{n-1}\nu(B_j) = \int\int\chi_{B_j}(x) dm_n(j) 
d\nu(x) \geq \int_X \theta d\nu(x)\geq\theta / 2 >0. 
\end{equation*}
This completes the proof of the lemma.
\end{proof}

A direct consequence of this lemma is 

\begin{corollary}The sequence $\lambda^{-n} g_n$ is uniformly bounded.
\end{corollary}

Next, we control the Holder norm of the sequence $g_n$:

\begin{lemma}There exists a constant $K_2>0$ such that 
\begin{equation*}
\left|\frac{g_n(x_1)}{g_n(x_2)}-1\right|\leq K_2 d(x_1,x_2)^{\alpha}
\end{equation*}
for any $x_1,x_2\in f(R_i)$ with $R_i\in\SR$. 
\end{lemma}

\begin{proof}Let $R^n$ be an atom of $\SR_h^n$ with $f^{n}(R^n)=R_i$. Then,
$x_1$ and $x_2$ have unique $f^n$ pre-images $y_1$ and $y_2$ inside $R^n$. By
corollary~\ref{c.3.7}, there exists a constant $A>0$ such that 
\begin{equation*}
|S_n\phi(y_1)-S_n\phi(y_2)|\leq A d(x_1,x_2)^{\alpha}.
\end{equation*}
Then, 
\begin{equation*}
e^{-A d(x_1,x_2)^{\alpha}}\leq \frac{g_n(x_1)}{g_n(x_2)}=\frac{\sum\limits_{R^n}
e^{S_n\phi(y_1)}}{\sum\limits_{R^n} e^{S_n\phi(y_2)}}\leq e^{A
d(x_1,x_2)^{\alpha}}
\end{equation*}
Now since $|e^{\pm A d(x_1,x_2)^{\alpha}}-1|\leq K_2 d(x_1,x_2)^{\alpha}$ if
$K_2$ is large depending on the constant $A$ and the diameter of $M$, we proved
the lemma.
\end{proof}

\begin{corollary}The sequence $\lambda^{-n} g_n$ is equicontinuous.
\end{corollary}

\begin{proof}This is easy from the previous lemma which says that the sequence
$\lambda^{-n} g_n$ has uniformly bounded H\"older norm.
\end{proof}

Applying the theorem of Arzela-Ascoli, we conclude that 
\begin{equation*}
\frac{1}{n}\sum\limits_{j=0}^{n-1} \lambda^{-i} g_i
\end{equation*}
has some subsequence converging uniformly on $M$ to an H\"older function $h$. We
claim that the function $h$ is an eigenfunction of the transfer operator (with
eigenvalue $\lambda$):

\begin{lemma}$\SL_{\phi}h=\lambda h$ and $h$ is bounded away from zero and
infinity.
\end{lemma}

\begin{proof}Since $\lambda^{-n} g_n\leq K_1$ and
$\frac{1}{n}\sum\limits_{j=0}^{n-1}\lambda^{-i} g_i\geq K_1^{-1}$, we conclude
that $K_1^{-1}\leq h\leq K_1$. 

To show that $\SL_{\phi}h=\lambda h$, note that 
\begin{equation*}
\SL_{\phi}h=\lim\limits_k\frac{1}{n_k}\sum\limits_{j=0}^{n_k-1}\lambda^{-j} 
(\SL_{\phi}g_j - g_{j+1}) + \frac{\lambda}{n_k}\sum\limits_{j=0}^{n_k-1}
\lambda^{-j}g_j - \frac{\lambda}{n_k} + \frac{\lambda^{-n_k}g_{n_k}}{n_k}
\end{equation*} 
Observe that $\frac{\lambda^{-n_k}g_{n_k}}{n_k}$ goes to zero. Thus, the lemma
is proved if we are able to obtain that 
\begin{equation*}
\frac{1}{n}\sum\limits_{j=0}^{n-1}\lambda^{-j} 
(\SL_{\phi}g_j - g_{j+1})
\end{equation*}
converges to zero. To begin with, denote by $E_{n+1}$ the set of cylinders 
\begin{equation*}
E_{n+1}:=\{[i_0,\dots,i_{n}]\in\SR^{n+1}; [i_0,\dots,i_{n}]\notin\SR_h^{n+1}
\textrm{ and } [i_1,\dots,i_n]\in\SR_h^n\}.
\end{equation*} 
Clearly $\# E_{n+1}\leq \#\{[i_0,\dots,i_{n}]\in\SR^{n+1}; 
\frac{\#\{j:i_j\leq q\}}{n+1}\geq\gamma\}\leq c_0^{n+1}$ for all large $n$.
Hence, 
\begin{equation*}
\|\lambda^{-j}(\SL_{\phi} g_j - g_{j+1})\|_{\infty}\leq
\lambda^{-i}\sum\limits_{R\in E_{j+1}} e^{S_j\phi(R)}\leq c_0^{j+1}\lambda^{-j}
e^{j\max\phi}
\end{equation*}
for any large $j$. Since $\lambda>e^{\max\phi+c_0}$, 
\begin{equation*}
\|\lambda^{-j}(\SL_{\phi} g_j - g_{j+1})\|_{\infty}\leq
\eta^i,
\end{equation*}
where $\eta<1$ and $j$ large. This completes the proof.
\end{proof} 

Next, we prove that the $f$-invariant measure $\mu = h\nu$ is mixing: 

\begin{lemma}\label{l.exact-gap}$\mu$ is exact.
\end{lemma}

\begin{proof}We begin with the following abstract claim: 
\begin{claim}\label{claim1}
Let $X$ be a compact metric space, $\mu$ a measure and 
$\SP=\{P_1,\dots,P_l\}$ a measurable partition of $X$. Consider a sequence of
collections $C_m$ of pairwise disjoint sets such that $\textrm{diam}(C_m):=
\max\limits_{C\in C_m} \textrm{diam}(C)\to 0$ and 
$\mu(X -\bigcup\limits_{C\in C_m} C)\to 0$ for all $m\geq m_0$ as $m\to\infty$.
Then, there exists $\{E_1^{(m)},\dots,E_l^{(m)}\}$ satisfying 
\begin{itemize}
\item Each $E_i^{(m)}$ is the union of atoms of $C_m$;
\item $\lim\limits_{m\to\infty} \mu(E_i^{m}\Delta P_i) = 0$ for all $1\leq i\leq
l$.
\end{itemize}  
\end{claim}
\begin{proof}[Proof of the claim~\ref{claim1}]Fix $\varepsilon>0$ and let
$K_1,\dots,K_l$ compact sets with $K_i\subset P_i$ and
$\mu(P_i-K_i)<\varepsilon$. Put $\delta:=\inf\limits_{i\neq j} d(K_i,K_j)$ and
consider $m$ large such that $\textrm{diam}(C_m)<\delta / 2$. We separate the
elements $C\in C_m$ into some collections $\SE_i^{(m)}$ characterized by
$C\in\SE_i^{(m)}$ if $C\cap K_i\neq\emptyset$. Note that any $C\in C_m$
intersects at most one of the compact sets $K_i$ (and if it does not intersect
one of them, we include $C$ arbitrarily into some $\SE_i^{(m)}$). Define
$E_i^{(m)}:=\bigcup\limits_{C\in\SE_i^{(m)}} C$. Then, 
\begin{equation*}
\begin{split}
&\mu(E_i^{(m)}\Delta P_i) = \mu(P_i - E_i^{(m)}) + \mu(E_i^{(m)} - P_i)\leq \\ 
&\leq \mu(P_i-K_i)+\mu(X-\bigcup\limits_{C\in C_m} C) + 
\mu(X-\bigcup\limits_{j=1}^l K_j)\leq (l+2)\varepsilon,
\end{split}
\end{equation*}
if $m$ is large enough so that $\mu(X-\bigcup\limits_{C\in C_m} C) < 
\varepsilon$. This proves the claim.
\end{proof}
Now, we are ready to show that $\mu$ is exact. If $\mu$ is not exact, there is
$A\in\bigcap\limits_{n\geq 0}f^{-n}(\SF)$ (where $\SF$ is the Borel
sigma-algebra) such that $\mu(A)>0$ and $\mu(M-A)>0$.
Using the claim~\ref{claim1}, we have that for all $\varepsilon>0$, there exists
$N(\varepsilon)$ such that 
\begin{equation*}
\frac{\mu(A\cap R^n)}{\mu(R^n)}\geq 1-\varepsilon,
\end{equation*}
for some $[i_0,\dots,i_{n-1}]=R^n\in\SR^n_h$ whenever $n\geq N(\varepsilon)$ 
(since the hyperbolic n-cylinders forms a collection of sets fitting 
the hypotheisis of the claim~\ref{claim1} by the corollary~\ref{c.3.6} 
and lemma~\ref{l.6.2}).

By definition, $A=f^{-n}(A_n)$ where $A_n\in\SF$. Because the distortion on
hyperbolic cylinders is bounded (see the corollary~\ref{c.3.8}), it holds 
\begin{equation*}
\frac{\mu(A_n\cap R_{i_0})}{\mu(R_{i_0})}\geq 1-\varepsilon K.
\end{equation*}
Similarly, 
\begin{equation*}
\frac{\mu(A_n^c\cap R_{j_0})}{\mu(R_{j_0})}\geq 1-\varepsilon K.
\end{equation*}
Since the hypothesis (H1) ensures that $\SR$ is a transitive partition, there is
a $N$ such that $f^N(R_i)=M$ for all $1\leq i\leq p+q$. In particular, the
image $Q_i$ of the rectangle $R_i$ under some inverse branch of $f^N$ is 
contained in $R_1$ for all $i$. Take $\delta<1/2K$ and $\varepsilon$ such that 
\begin{equation*}
\varepsilon\cdot \frac{K\sup_i \mu(R_i)}{c}<\delta,
\end{equation*}
where $c=\min_i\mu(Q_i)$. So, 
\begin{equation*}
\frac{\mu(A_n\cap Q_{i_0})}{\mu(Q_{i_0})}, 
\frac{\mu(A_n^c\cap Q_{j_0})}{\mu(Q_{j_0})}>1-\delta.
\end{equation*}
Since $A_n=f^{-N}(B_N)$ and $A_n^c=f^{-N}(B_N^c)$, the bounded distortion result
of the corollary~\ref{c.3.8} implies 
\begin{equation*}
\frac{\mu(B_N\cap R_{i_0})}{\mu(R_{i_0})}, 
\frac{\mu(B_N^c\cap R_{j_0})}{\mu(R_{j_0})} > 1-\delta K.
\end{equation*}
Summing up these estimates, we get $1>2-2\delta K$, a contradiction.
\end{proof}

To proceed further, we introduce the compact operators $T_n$ given by 
\begin{equation*}
T_n  = \SL^n_{\phi}\circ \pi_n ,
\end{equation*}
where $\pi_n g(x):=\frac{1}{\nu(R^n(x))}\int_{R^n(x)} g d\nu$. Note that $\pi_n$
is well-defined since $\nu(R^n)>0$ for any $R^n\in\SR^n$, \footnote{This follows
directly from $J_{\nu}f=\lambda e^{-\phi}\leq C:=k e^{\max\phi}$ and 
$0<\nu(f(R_{i_{n-1}})) = \nu(f^n(R^n)) = \int_{R^n} J_{\nu}f^n\leq 
C^n \nu(R^n)$.} and $T_n$ 
has finite dimensional rank since $\SR^n$ is a finite partition. 

In the sequel, we show the following result about the approximation of
$\SL_{\phi}$ by compact operators: 

\begin{lemma}\label{l.approx-gap}There are constants $C>0$ and $0<\theta<1$ such
that 
\begin{equation*}
\frac{1}{\lambda^n}\|\SL_{\phi}^n g - T_n g\|_{\infty}\leq C\theta^n
|||g|||_{\alpha}
\end{equation*}
and 
\begin{equation*}
\frac{1}{\lambda^n}|||\SL_{\phi}^n g - T_n g|||_{\alpha}\leq C\theta^n
|||g|||_{\alpha}.
\end{equation*}
\end{lemma}

\begin{proof} Given $g$ and $R^n_i\in\SR^n$, we fix $z_i\in R^n_i$ with 
$\pi_n g(z)=g(z_i)$ for all $z\in R^n_i$. Also, given $x\in M$, we denote by 
$y_i$ the pre-image of $x$ under $f^n$ such that $y_i\in R^n_i$. We have 
\begin{equation*}
\begin{split}
&\frac{1}{\lambda^n}|\SL_{\phi}^n g(x) - T_n g(x)|\leq \frac{1}{\lambda^n}
\sum\limits_{R^n_i\in\SR^n} 
e^{S_n\phi(y_i)}|g(y_i) - g(z_i)| \leq \\ 
&\frac{e^{n\max\phi}}{\lambda^n} \left(
\sum\limits_{R^n_i\in A(n)} |g(y_i)-g(z_i)| +   
\sum\limits_{R_i^n\notin A(n)} |g(y_i)-g(z_i)|\right), 
\end{split}
\end{equation*}
where $A(n)$ is the collection of $n$-cylinders $R^n=[i_0,\dots, i_{n-1}]$ with 
\begin{equation*}
\frac{1}{n}\#\{0\leq j < n; \, i_j\leq q\}> \gamma.
\end{equation*}
It follows that 
\begin{equation*}
\begin{split}
\frac{1}{\lambda^n}|\SL_{\phi}^n g(x) - T_n g(x)|&\leq  
\frac{\left(e^{(c_0+\max\phi) n} (1+\delta_0)^n + k^n
e^{(-c+\max\phi)n}\right)}{\lambda^n} |||g|||_{\alpha}\leq \\ 
&\leq C\theta^n |||g|||_{\alpha},
\end{split}
\end{equation*}
for some $0<\theta<1$ (if $\delta_0$ and $\max\phi$ are sufficiently small).
Similarly, 
\begin{equation*}
\begin{split}
&\frac{1}{\lambda^n}|||\SL_{\phi}^n g - T_n g|||_{\alpha}\leq \\  
&\frac{1}{\lambda^n}\left((1+\delta_0)^n e^{c_0 n} +
e^{-cn} k^n\right)
\left(e^{n\max\phi} + 
2e^{n\max\phi} n\max\phi\right)|||g|||_{\alpha} \leq \\ 
&C\theta^n |||g|||_{\alpha}.
\end{split}
\end{equation*}
This completes the proof of the lemma.
\end{proof}

Finally, the spectral gap property of the normalized transfer operator
$\widetilde{\SL_{\phi}}:=\frac{1}{\lambda}\SL_{\phi}$ (i.e.,
$\textrm{spec}\widetilde{\SL_{\phi}}=\{1\}\cup\Sigma_0$ where $1$ is a simple
eigenvalue and $\Sigma_0$ is contained in a disc of radius $\theta<1$) is a
direct corollary of the two previous lemmas. Indeed, the
lemma~\ref{l.approx-gap} implies that 
$\lambda^{-n}\|\SL_{\phi}^n - T_n\|_{\alpha}\leq C\theta^n$ for some $\theta<1$,
where $T_n$ are compact operators. Hence, the spectral radius of $\SL_{\phi}$ is
$\theta<1$ at most. So, we can write $\textrm{spec}(\widetilde{\SL}) =
\{1=\beta_0,\dots,\beta_m\}\cup \Sigma_0$ with $\Sigma_0$ contained in a disk of
radius $\theta<1$ and $\beta_i$ are eigenvalues with norm $1$. Now, the
lemma~\ref{l.exact-gap} allow us to prove that $m=0$, i.e., there are no
eigenvalues with norm $1$ other than $\beta_0=1$. In fact, if there exists some
non-zero H\"older function $v_i$ with $\widetilde{\SL_{\phi}}v_i =
\beta_i v_i$ and $|\beta_i|=1$, then the exactness property of $\mu$ (see
lemma~\ref{l.exact-gap}) implies 
\begin{equation*}
\begin{split}
\int u\widetilde{\SL_{\phi}}^n v_i d\nu = \int (u\circ f^n) v_i d\nu = 
\int (u\circ f^n) (v_i/h) d\mu \to \\ 
\to \int u d\mu \int (v_i/h) d\mu = \int u \left(h\int v_i d\nu\right) d\nu 
\end{split}
\end{equation*} 
so that $\beta_i^n v_i = \widetilde{\SL_{\phi}}^n v_i$ converges to $h\int v_i
d\nu$. Hence, $\beta_i=1$ and $v_i = h\int v_i d\nu$. This proves that $1$ is
the only eigenvalue of $\widetilde{\SL_{\phi}}$ with norm $1$ and its eigenspace
is one-dimensional. Furthermore, its algebraic multiplicity is one since,
otherwise there is $u_0$ with $\widetilde{\SL_{\phi}}^n u_0 = h + n u_0$, a
contraction with the boundedness of $\|\widetilde{\SL_{\phi}}^n\|_{\alpha}$.
Note that the spectral splitting of $\widetilde{\SL_{\phi}}$ is $\mathbb{R} h
\oplus X_0$, where $X_0:=\{u: \int u d\nu =0\}$. 

To close this subsection, we show that the spectral gap property implies that
$\mu$ is exponentially mixing:

\begin{corollary}$\mu$ has exponential decay of correlations.
\end{corollary}

\begin{proof}Clearly, 
\begin{equation*}
\begin{split}
\int (u\circ f^n) v d\mu - \int u d\mu \int v d\mu = \int u
\widetilde{\SL_{\phi}}^n (\pi_0(v\cdot h)) d\nu,
\end{split}
\end{equation*}
where $\pi_0(\varphi):=\varphi - \int\varphi d\nu$ is the spectral projection
onto $X_0=\{u: \int u d\nu = 0\}$. Since 
\begin{equation*}
\begin{split}
&|\int u\widetilde{\SL_{\phi}}^n (\pi_0(v\cdot h)) d\nu|\leq
\|\widetilde{\SL_{\phi}}(\pi_0(v\cdot h))\|_{\infty} \int |u| d\nu \\ 
&\leq \|\widetilde{\SL_{\phi}}(\pi_0(v\cdot h))\|_{\alpha} \int |u| d\nu \leq
C\theta^n \|\pi_0(vh)\|_{\alpha}\int |u| d\nu \leq \\ 
&C'\int |u| d\nu \|v\|_{\alpha} \theta^n = K(u,v) \theta^n. 
\end{split}
\end{equation*}
\end{proof}

\subsection{Stochastic Stability}

The objective of this subsection is to prove the stability of the equilibrium
measures under small random perturbations of the maps and potentials considered
here. By small random pertubations of a fixed map $f_0$ we mean that a 
metric space $\Omega$, a transformation $\omega\in\Omega\mapsto f_{\omega}$ with
$f_{\omega}\to f_{\omega_0}=f_0$ in the $C^{1+\alpha}$ topology as 
$\omega\to\omega_0$ for some
$\omega_0\in\Omega$ and a family of probabilities $\theta_{\varepsilon}$ on
$\Omega$ such that $\textrm{supp}(\theta_{\varepsilon})\to\{\omega_0\}$ as
$\varepsilon\to 0$ are given. In a similar fashion, we can define random
perturbations of the potential $\phi$. 

To handle this problem, we mainly use the works of Baladi and
Young~\cite{B},~\cite{BY}. Namely, it follows from these papers that the
stochastic stability of the equilibrium states for the transformations $f_0$ 
and potentials $\phi$ in our hypothesis is guaranteed if we can prove the
fundamental lemma: 
 
\begin{lemma}
If $\delta_0$ in hypothesis (H1) is sufficiently small, 
there are constants $C>0$, $0<\theta<1$ and, for each $n$, some
$\varepsilon(n)>0$ such that for any $0<\varepsilon<\varepsilon(n)$, 
\begin{equation*}
\frac{1}{\lambda^n}|||\SL_{\varepsilon}^n g|||_{\alpha}\leq C\theta^n
|||g|||_{\alpha} + C\|g\|_{\infty}.
\end{equation*} 
and 
\begin{equation*}
\frac{1}{\lambda^n}\|\SL_{\varepsilon}^n g - \SL^n g\|_{\infty}\leq C\theta^n 
(\|g\|_{\infty}+ |||g|||_{\alpha}).
\end{equation*}
Here $\SL=\SL_{\phi}$ is the transfer operator of $(f_0,\phi)$ and 
$\SL_{\varepsilon}$ is the random version of the transfer operator when the
random noise has level $\varepsilon$.\footnote{See~\cite{B} for the precise 
definitions.}
\end{lemma}

\begin{proof}The first part of the lemma follows from a simple modification of
the arguments in the lemma~\ref{l.LY} using that the random noise is small.
Also, the second part is quite similar. Indeed, if the random noise is small, 
for every $x\in M$ and $y$ with $f^n(y)=x$, there exists a unique 
$y_{\varepsilon}$
close to $y$ such that $f_{\varepsilon}^n(y_{\varepsilon})=x$. Hence, 
\begin{equation*}
\begin{split}
\frac{1}{\lambda^n}\|\SL_{\varepsilon}^n g(x) - \SL^n g(x)\|_{\infty}&\leq
\frac{1}{\lambda^n}\sum\limits_{i} 
|e^{S_n\phi_{\varepsilon}(y^{(i)}_{\varepsilon})}g(y^{(i)}_{\varepsilon}) -
e^{S_n\phi(y^{(i)})}g(y^{(i)})| \\ 
&\leq C\theta^n\cdot |||g|||_{\alpha}+C\xi^n\cdot \|g\|_{\infty}, 
\end{split}
\end{equation*}
where $\theta,\xi\in (0,1)$ satisfies 
$$1>\theta^n> 
\frac{1}{\lambda^n}e^{(c_0+\max\phi_{\varepsilon})n}(1+\delta_0(\varepsilon))^n
+ \frac{1}{\lambda^n}e^{n\max\phi_{\varepsilon}}e^{-4cn} k^n$$ and 
$$1>\xi^n>\frac{2}{\lambda^n}e^{n\max\phi}n\max\phi$$ with 
$\sigma_1(\varepsilon)$, $\delta_0(\varepsilon)$ such that  
$\sigma_1(\varepsilon)\to\sigma_1$, $\delta_0(\varepsilon)\to\delta_0$
as $\varepsilon\to 0$. 
Of course, the smallness of the random
noise $(f_{\varepsilon},\phi_{\varepsilon})$ was used to get this estimate. This
finishes the proof since our hypothesis (besides the facts 
$f_{\varepsilon}\to f$
and $\phi_{\varepsilon}\to\phi$ as $\varepsilon\to 0$) ensures the existence of
$\theta,\xi<1$ as above. 
\end{proof}

\begin{remark}In connection with the results of Arbieto, Matheus and 
Oliveira, the
stochastic stability theorem proved here implies that the quenched equilibrium
states of~\cite{AMO} are unique and exponentially mixing. 
\end{remark}

\subsection{Large Deviations} 

Given a transformation $f:M\to M$ and an invariant measure $\mu$, we say that
$(f,\mu)$ has the \emph{large deviations property} if for any continuous 
observable $\varphi$ and $\rho>0$, 
\begin{equation*}
\limsup\limits_{n\to\infty}\frac{1}{n}\log\mu\left(\left\{x\in M:
\frac{1}{n}S_n\varphi(x)\notin \left(\int\varphi d\mu-\rho,
\int\varphi d\mu+\rho\right)\right\}\right) <0.
\end{equation*} 

The goal of this subsection is to prove the following result: 

\begin{theorem}The equilibrium states $\mu$ of the $(f,\phi)$ in both 
theorems~\ref{t.A} and~\ref{t.B} have the large deviation property:
\begin{equation*}
\begin{split}
&\limsup\limits_{n\to\infty}\frac{1}{n}
\log\mu(\{x:\frac{1}{n}S_n\varphi\notin (\int\varphi d\mu -\rho, 
\int\varphi d\mu + \rho)\})\leq \\ 
&\sup\{h_{\eta}(f)+\int\phi d\eta - P: \eta
\text{ is } f-\text{invariant} \text{ and } |\int\varphi d\eta - \int\varphi
d\mu|\geq\rho\} \\ 
&<0.
\end{split}
\end{equation*}
\end{theorem}

To aleviate the notation during the proof of this theorem, we set 
\begin{equation*}
B_n:=\{x\in M: \frac{1}{n}S_n\varphi(x)\geq \rho\}
\end{equation*} 
Since $\varphi$ is continuous, we can find
$E\subset U\subset B_n$ such that $U$ is open and
$\nu(B_n-E)<\frac{1}{2n}\nu(B_n)$. Using the lemma~\ref{l.6.2}, there is a
family $\mathcal{F}_n=\mathcal{E}_1\cup\dots\cup\mathcal{E}_k$ of hyperbolic
cylinders contained in $U$ verifying 
\begin{equation*}
\nu(E\triangle \bigcup\limits_{F\in\mathcal{F}_n} F) \leq
\left(1-\frac{\theta}{4}\right)^k<\frac{1}{2n}\nu(B_n).
\end{equation*}

Note that each $P\in\mathcal{E}_i$ there is $x\in M$ and $h_i$ an hyperbolic
time for $x$ such that $P=R^{h_i}(x)$, where $i=1,\dots,k$ and
$n<h_1<h_2<\dots<h_k$. Define $\mathcal{C}_n$ the set of all such pairs
$(x,h_i)$, one for each element of $\mathcal{F}_n$. 

Let 
\begin{equation*}
\sigma_n = \frac{1}{Z_n}\sum\limits_{(x,l)\in\mathcal{C}_n} e^{S_l\xi(x)}
\cdot \delta_x, \text{ where } 
Z_n=\sum\limits_{(x,l)\in\mathcal{C}_n} e^{S_l\xi(x)} \text{ and } 
\xi=\phi-P.
\end{equation*}
Using the lemma 9.9 in~\cite{W}, we obtain 
\begin{equation*}
H_{\sigma_n}(\bigvee\limits_{i=0}^{h_k - 1} f^{-i}(\mathcal{R})) + \int
S_{l(x)}\xi(x) d\sigma_n(x) = \log\sum\limits_{(x,l)\in\mathcal{C}_n}
e^{S_l\xi(x)}, 
\end{equation*}
where $l(x)$ is the unique integer $l$ such that $(x,l)\in\mathcal{C}_n$. Since
$\max\phi < P$, it follows that $S_{l(x)-n}\xi(x)<0$ whenever $l(x)>n$. 
Therefore, 
\begin{equation*}
H_{\sigma_n}(\bigvee\limits_{i=0}^{h_k - 1} f^{-i}(\mathcal{R})) + \int
S_{n}\xi(x) d\sigma_n(x) \geq \log\sum\limits_{(x,l)\in\mathcal{C}_n}
e^{S_l\xi(x)}
\end{equation*}
Putting $\eta_n=\frac{1}{n}\sum\limits_{i=0}^n f^i_* \sigma_n$ and $\eta$ any 
accumulation point of $\eta_n$, we can repeat the arguments in p. 220 
of~\cite{W} to derive that 
\begin{equation*}
\limsup\limits_{n\to\infty}\frac{1}{n}\log Z_n \leq h_{\eta}(f,\mathcal{R}) +
\int\xi d\eta\leq h_{\eta}(f)+\int\phi d\eta - P.
\end{equation*}

Observe that $\eta_n$ is a convex combination of measures of the form
$\frac{1}{n}\sum\limits_{i=0}^{n-1}\delta_{f^i(x)}$. In particular, 
\begin{equation*}
\int\varphi d\eta_n = \frac{1}{Z_n}\sum\limits_{(x,l)\in\mathcal{C}_n}
e^{S_l\xi(x)}\cdot\frac{1}{n}\sum\limits_{i=0}^{n-1}\varphi(f^i(x))\geq
\frac{1}{Z_n}\sum\limits_{(x,l)\in\mathcal{C}_n}
e^{S_l\xi(x)}\cdot \rho = \rho.
\end{equation*} 
So, $\int\varphi d\eta\geq\rho$ since $\varphi$ is continuous. 

Also, the definition of $\mathcal{C}_n$ and the non-lacunary weak Gibbs property
of $\nu$ (see the lemma~\ref{l.3.12}) implies 
\begin{equation*}
\begin{split}
\nu(B_n)&\leq \frac{n}{n-1}\nu(\bigcup\limits_{F\in\mathcal{F}_n})\leq
\frac{n}{n-1}\sum\limits_{(x,l)\in\mathcal{C}_n} \nu(R^l(x)) \\ 
&\leq
\frac{n}{n-1}\sum\limits_{(x,l)\in\mathcal{C}_n}Ke^{S_l\xi(x)} = \frac{Kn}{n-1}
Z_n.
\end{split}
\end{equation*}

Hence, we constructed an invariant measure $\eta$ with $\int \varphi d\eta\geq
\rho$ and 
\begin{equation*}
\limsup\limits_{n\to\infty}\frac{1}{n}\log\nu(B_n)\leq
\limsup\limits_{n\to\infty}\frac{1}{n}\log Z_n\leq h_{\eta}(f)+\int\phi d\eta -
P.
\end{equation*}

This estimate allows us to prove the large deviations property for our
equilibrium state $\mu$. Indeed, since $\mu= h\nu$ with $h$ (H\"older)
continuous, we showed that for any continuous function 
$\varphi$ and any fixed number $\rho>0$, 
\begin{equation*}
\begin{split}
&\limsup\limits_{n\to\infty}\frac{1}{n}
\log\mu(\{x:\frac{1}{n}S_n\varphi\notin (\int\varphi d\mu -\rho, 
\int\varphi d\mu + \rho)\})\leq \\ 
&\sup\{h_{\eta}(f)+\int\phi d\eta - P: \eta
\text{ is } f-\text{invariant} \text{ and } |\int\varphi d\eta - \int\varphi
d\mu|\geq\rho\}.
\end{split}
\end{equation*}

Thus, our task is to show that the right-hand side of this estimate is strictly
negative. Note that we have shown in the proof of proposition~\ref{l.6.5} 
that $h_{\eta}(f)+\int\phi d\eta-P\leq 0$ for any $f$-invariant measure $\eta$.
If the supremum in the right side is zero, there exists a sequence of
$f$-invariant measures $\eta_n$ with $|\int\varphi d\eta_n-\int\varphi
d\mu|\geq\rho$ and $h_{\eta_n}(f)+\int\phi d\eta-P\to 0$. The proof of the
Rokhlin's formula in the appendix III (see corolaries~\ref{c.rokhlinA}
and~\ref{c.rokhlinB} below) says that the $\eta_n$ are expanding measures for
sufficiently large $n$, so that the lemma~\ref{l.appendixIII} can be used to give
that any partition with small diameter is generating for any such $\eta_n$.
Taking any accumulation point $\eta$ of the sequence $\eta_n$, then we have
$h_{\eta_n}(f)\to h_{\eta}(f)$. Therefore, since the potential $\phi$ is
continuous, 
\begin{equation*}
h_{\eta}(f)+\int\phi d\eta-P = \lim\limits_{n\to\infty} 
h_{\eta_n}(f)+\int\phi d\eta_n - P = 0.
\end{equation*} 
In particular, $\eta$ is an equilibrium measure. By uniqueness, the measures 
$\eta$ and $\mu$ coincides. But this a contradiction since $\int\varphi d\eta =
\lim\int\varphi d\eta_n\notin (\int\varphi d\mu -\rho, \int\varphi d\mu +
\rho)$. This completes the proof of the large deviations property of $\mu$.

\bigskip 

\textbf{Acknowledgements.} The authors are thankful to Krerley Oliveira and
Marcelo Viana for sharing their ideas about the thermodynamical formalism of
non-uniformly expanding maps contained in the papers~\cite{OV1},~\cite{OV2}. 
Also, we are grateful to F. Ledrappier for the suggestion of a 
Ruelle-Perron-Frobenius transfer operator approach to the study of the speed of
mixing of the equilibrium states considered here and C. Liverani for the
discussions around the explict expression of the family of inavriant cones
$\Lambda_L$ used in this paper. Furthermore, we are indebted 
to Viviane Baladi for the invitation to attend the ``Trimester Time at Work''
(May-June 2005) at the Institute Henri Poincar\'e where this project started. 
Finally, the authors are thankful to IMPA and its staff for the fine scientific 
ambient. 

\section{Appendix I: a combinatorial lemma}

Consider the set 
\begin{equation*}
I_{\gamma,n}:=\left\{(i_0,\dots,i_{n-1})\in \{1,\dots,q,\dots,q+p\}^{n}; \; 
\#\{0\leq j<n; \; i_j\leq q\}>\gamma n\right\}
\end{equation*}

By definition,
\begin{equation*}
\#I_{\gamma,n}\leq \sum\limits_{r\geq\gamma n}\binom{n}{r} p^{n-r} q^r.
\end{equation*}

The general plan of the proof of lemma~\ref{l.2.7} is: our goal is to 
estimate the rate of exponential growth of 
$I_{\gamma,n}$ (when $n\to\infty$) for $\gamma$ close to $1$. A natural way to 
do this is to use the Stirling's formula to compare the binomial numbers with
exponentials functions, in order to get some good asymptotics.

Now, let us see how this strategy works. A direct consequence of Stirling's 
formula is the existence of a 
universal constant $B>0$ such that:
\begin{equation*}
\binom{n}{r} \leq
B\left((1+\frac{1}{k})(1+k)^{1/k}\right)^r
\leq B\left( (1+k)^{1/k}\cdot (1+\frac{1}{k})\right)^n.
\end{equation*}
for every $r\geq \frac{k}{k+1} n$. So, if $\gamma\geq \frac{k}{k+1}$, then
\begin{equation*}
c_{\gamma}\leq \log\left((1+\frac{1}{k})(1+k)^{1/k}p^{1/(k+1)}q\right).
\end{equation*}

Since the right-hand side of this inequality goes to $\log q$ when 
$k\to\infty$, this completes the proof of lemma~\ref{l.2.7}.

\section{Appendix II: positive operators and cones}\label{a.cones}

This appendix presents some results of the theory of projective metrics on cones
and positive operators (due to
Garrett Birkhoff) used in the subsection~\ref{s.s.cones}.\footnote{For sake of
simplicity, we consider cones and positive operators
on Banach spaces only (which is precisely the context of 
subsection~\ref{s.s.cones}); 
however, the arguments can be extended to the general case of vector spaces 
with a partial ordering satisfying some mild conditions. See~\cite{L} for
more details.}

Let $\SB$ be a Banach space with the topology of the norm. A subset
$\Lambda\subset \SB-\{0\}$ is a \emph{cone} if $r\cdot v\in\Lambda$ for all
$v\in\Lambda$ and $r\in\real^+$. The cone $\Lambda$ is \emph{closed} if
$\Lambda\cup\{0\}$ and $\Lambda$ is \emph{convex} if $v+w\in\Lambda$ for all 
$v,w\in\Lambda$. 

A convex cone $\Lambda$ with $\Lambda\cap (-\Lambda) = \emptyset$ determine a 
\emph{partial ordering} $\preceq$ on $\SB$: 
\begin{equation*}
w\preceq v \; \text{ iff } \; v-w\in\Lambda\cup\{0\}.
\end{equation*} 

In the sequel, our cones $\Lambda$ are assumed to be closed, convex and
$\Lambda\cap (-\Lambda)=\emptyset$.

Given a cone $\Lambda$ and two vectors $v,w\in\Lambda$, we define 
$\Psi (v,w)=\Psi_{\Lambda}(v,w)$ by 
\begin{equation*}
\Psi (v,w)=\log\frac{B_{\Lambda}(v,w)}{A_{\Lambda}(v,w)},
\end{equation*}
where $A_{\Lambda}(v,w)=\sup\{r\in\real^+: \; r\cdot v\preceq w\}$ and
$B_{\Lambda}(v,w)=\inf\{r\in\real^+: \; w\preceq r\cdot v\}$. The
(pseudo-)metric $\Psi$ is called the \emph{projective metric} of $\Lambda$ (or
\emph{$\Lambda$-metric} for brevity).\footnote{The justification of the name
``projective metric'' for the pseudo-metric $\Psi$ is: defining the equivalence
relation $v\sim w$ iff $w=r\cdot v$ for some $r\in\real^+$, then $\Psi$
induces a metric on the quotient $\Lambda /\sim$.}

A key result of this appendix is:

\begin{theorem}\label{t.Birkhoff}
Let $\Lambda_i$ be a closed convex cone (with $\Lambda_i\cap (-\Lambda_i) =
\emptyset$) in a Banach space $\SB_i$, for $i=1,2$. If $\SL:\SB_1\to\SB_2$ is a
linear operator such that $\SL(\Lambda_1)\subset\Lambda_2$, then 
\begin{equation*}
\Psi_{\Lambda_2} (\SL v,\SL w)\leq
\tanh\left(\frac{\text{diam}_{\Psi_{\Lambda_2}}(\SL\Lambda_1)}{4}\right)\cdot 
\Psi_{\Lambda_1} (v,w),
\end{equation*}
for any $v,w\in\Lambda$.
\end{theorem}

\begin{proof}Take $v,w\in \Lambda_1$. If $A_{\Lambda_1}(v,w)=0$ or
$B_{\Lambda_1}(v,w)=\infty$, the desired inequality is easy. On the other hand, 
if $A:=A_{\Lambda_1}(v,w)\neq 0$ and $B:=B_{\Lambda_1}(v,w)\neq\infty$, then 

\begin{equation*}
\Psi_{\Lambda_1} (v,w)=\log\frac{B}{A}
\end{equation*}
and $A\cdot v \preceq w$ and 
$w\preceq B\cdot v$. Observe that $A_{\Lambda_1}(v,w)$ and 
$B_{\Lambda_1}(v,w)$ are non-negative. In particular, if
$\text{diam}_{\Psi_{\Lambda_2}}(\SL\Lambda_1)=\infty$, the desired inequality
follows. If $\Delta:=\text{diam}_{\Psi_{\Lambda_2}}(\SL\Lambda_1)<\infty$, then 

\begin{equation*}
\Psi_{\Lambda_2} (\SL (w-Av),\SL (Bv-w))\leq \Delta.
\end{equation*}
Hence, there are $r,s\geq 0$ such that 

\begin{equation*}
r\cdot\SL (w-Af)\preceq \SL (Bv-w)
\end{equation*}
and

\begin{equation*}
\SL (Bv-w)\preceq s\cdot\SL (w-Av) 
\end{equation*}
with $\log\frac{s}{r}\leq\Delta$. Thus, 

\begin{equation*}
\SL w\preceq \frac{B+rA}{1+r}\cdot\SL v
\end{equation*} 
and

\begin{equation*}
\frac{sA+B}{1+s}\cdot\SL v\preceq \SL w
\end{equation*}
Therefore, 

\begin{equation*}
\begin{split}
\Psi_{\Lambda_2} (\SL v,\SL w) &\leq \log\frac{(B+rA)(1+s)}{(sA+B)(1+r)} =
\log\frac{r+e^{\Psi_{\Lambda_1}(v,w)}}{s+e^{\Psi_{\Lambda_1}(v,w)}} -
\log\frac{1+r}{1+s} \\
&=\int_0^{\Psi_{\Lambda_1}(v,w)}\frac{(s-r)e^t}{(r+e^t)(s+e^t)}dt \leq 
\Psi_{\Lambda_1} (v,w)\frac{1-\frac{r}{s}}{\left(1+\sqrt{\frac{r}{s}}\right)^2}
\\ 
&\leq \tanh\left(\frac{\Delta}{4}\right)\Psi_{\Lambda_1} (v,w).
\end{split}
\end{equation*}
This finishes the proof.
\end{proof}

\begin{lemma}\label{l.appendixII} 
Let $\Lambda$ be a closed convex cone (with $\Lambda\cap
(-\Lambda)=\emptyset$) in a Banach space $\SB$ endowed with two norms
$\|.\|_{(i)}$, $i=1,2$ (not necessarily equivalent), and assume that for all
$v,w\in\SB$,
\begin{equation*}
-v\preceq w\preceq v \quad \text{ implies } \quad \|w\|_{(i)}\leq \|v\|_{(i)}
\quad \text{ for } \quad i=1,2. 
\end{equation*}
Then, for all $v,w\in\Lambda$ with $\|v\|_{(1)}=\|w\|_{(1)}$, we have 
\begin{equation*}
\|v-w\|_{(2)}\leq (e^{\Psi_{\Lambda}(v,w)}-1)\|v\|_{(2)}.
\end{equation*}
\end{lemma}

\begin{proof} $\Psi (v,w) = \log\frac{B}{A}$, where $Av\preceq w\preceq Bv$. In
particular, $-w\preceq 0\preceq Av\preceq w$ and, hence, $A\|v\|_{(1)}\leq
\|w\|_{(1)}$. Since $\|v\|_{(1)}=\|w\|_{(1)}$, we have $A\leq 1$. Similarly, it
is not hard to obtain $B\geq 1$. Thus,

\begin{equation*}
-(B-A)v\preceq (A-1)v\preceq w-v\preceq (B-1)v\preceq (B-A)v
\end{equation*} 
As a consequence, we have  
$$\|w-v\|_{(2)}\leq (B-A)\|v\|_{(2)}\leq 
\frac{B-A}{A}\|v\|_{(2)} = (e^{\Psi(v,w)}-1)\|v\|_{(2)}.$$

This concludes the proof of the lemma.
\end{proof}

\section{Appendix III: Rokhlin's formula}\label{a.rokhlin}

This appendix has a lot of non-trivial intersection with the recent work of
Oliveira and Viana~\cite{OV1}. In particular, the proposition~\ref{p.rokhlin} 
(and the lemmas~\ref{l.rokhlin},~\ref{l.rokhlin2} used in its proof) and the
lemma~\ref{l.appendixIII} below were borrowed from~\cite{OV1}.

We start with an abstract criterion concerning Rokhlin's formula for certain
measures with \emph{generating partitions}.
  
Let $f:M\to M$ be a measurable
transformation and $\mu$ an invariant probability. Suppose that there exist a
finite or countable partition $\SP$ such that 

\begin{itemize} 
\item (a) $f$ is locally injective (i.e.,
$f$ is injective on every atom of $\SP$); 
\item (b) $\SP$ is $f$-generating with
respect to $\mu$ (i.e., $\text{ diam }(\SP^n(x))\to 0$ for $\mu$ almost every
$x$, where $\SP^n(x)$ is the atom of $x$ in the partition
$\SP^n:=\bigvee\limits_{j=0}^{n-1}f^{-j}(\SP)$).
\end{itemize}

\begin{proposition}\label{p.rokhlin}If $\mu$ is an invariant measure satisfying 
(a) and (b) above, then $\mu$ verifies Rokhlin's formula
\begin{equation*}
h_{\mu}(f)=\int \log J_{\mu}f.
\end{equation*} 
\end{proposition}

\begin{proof}Let $\SP_{\infty}=\bigvee\limits_{j=0}^{\infty}f^{-j}(\SP)$ and
$\SQ_{\infty}=\bigvee\limits_{j=1}^{\infty}f^{-j}(\SP)$. Note that
$\SQ_{\infty}(x)=f^{-1}(\SP_{\infty}(f(x)))$. On the other hand, $\SP$ is
generating means that $\SP_{\infty}(x)=\{x\}$, and so
\begin{equation*}
\SQ_{\infty}(x)=f^{-1}(f(x)),
\end{equation*}
for $\mu$-almost all $x$.

Denote by $\mathbb{E}_{\mu}(\varphi|\SN)$ the ($\mu$-) conditional expectation 
of a function $\varphi:M\to\real$ relative to a partition $\SN$, i.e., the 
essentially
unique $\SN$-measurable function $\mathbb{E}_{\mu}(\varphi|\SN)$ such that
\begin{equation*}
\int_B \mathbb{E}_{\mu}(\varphi|\SN) d\mu=\int_B\varphi d\mu,
\end{equation*}
for every $\SN$-measurable set $B$.

At this point we prove the following lemma:

\begin{lemma}\label{l.rokhlin}
$\mathbb{E}_{\mu}(\varphi|\SQ_{\infty})=\sum\limits_{y\in\SQ_{\infty}(x)}
\frac{1}{J_{\mu}f(y)}\varphi(y)$ for $\mu$-a.e. $x$.
\end{lemma}

\begin{proof}[Proof of the lemma~\ref{l.rokhlin}] It is clear that the
right-hand side is $\SQ_{\infty}$-measurable. Let $B$ be any
$\SQ_{\infty}$-measurable set, that is, $B$ consists of entire atoms of
$\SQ_{\infty}$. From the definitions, there exists a measurable set $C$ such that
$B=f^{-1}(C)$. Hence, the $f$-invariance of $\mu$ implies
\begin{equation*}
\begin{split}
\int_B \sum\limits_{y\in\SQ_{\infty}(x)}\frac{1}{J_{\mu}f(y)}\varphi(y) d\mu(x)
&= \int_C\sum\limits_{y\in f^{-1}(z)}\frac{1}{J_{\mu}f(y)}\varphi(y) d\mu(z) \\ 
=\sum\limits_{A\in\SP}\int_{C_A}\frac{1}{J_{\mu}f(y_A)}\varphi(y_A) d\mu(z),
\end{split}
\end{equation*}
where $C_A = f(B\cap A)$ and $y_A = (f|A)^{-1}(z)$. Since every $f|A$ is
injective, we can use the definition of Jacobian to rewrite the later expression
as
\begin{equation*}
\sum\limits_{A\in\SP}\int_{B\cap A}\varphi(y) d\mu(y) = \int_B \varphi d\mu.
\end{equation*}
This proves the lemma.
\end{proof}

Recall that the conditional entropy of the partition $\SP$ with respect to the
partition $\SN$ is defined as 
\begin{equation*}
H_{\mu}(\SP|\SN):= \int\sum\limits_{A\in\SP}
-\mathbb{E}_{\mu}(\chi_A|\SN)\log\mathbb{E}_{\mu}(\chi_A|\SN) d\mu.
\end{equation*}
See~\cite[definition 4.8]{W}. Using the previous lemma, we can calculate the
conditional entropy of $\SP$ with respect to $\SQ_{\infty}$ in terms of the
Jacobian of $\mu$ as follows:
\begin{lemma}\label{l.rokhlin2}
$$H_{\mu}(\SP|\SQ_{\infty})=\int\log J_{\mu}f d\mu.$$
\end{lemma}
\begin{proof}[Proof of lemma~\ref{l.rokhlin2}]The lemma~\ref{l.rokhlin} says
that
\begin{equation*}
\mathbb{E}_{\mu} (\chi_A|\SQ_{\infty}) = \psi_A\circ f, \quad \text{ where }
\quad \psi_A(z)=\sum\limits_{y\in f^{-1}(z)} \frac{1}{J_{\mu} f(y)}\chi_A(y).
\end{equation*}
Observe that if $z\in f(A)$, then $\psi_A(z)=1/J_{\mu}f(y_A)$, where
$y_A=(f|A)^{-1}(z)$ and if $z\notin f(A)$, then $\psi_A(z)=0$. In particular,
\begin{equation*}
\begin{split}
H_{\mu}(\SP|\SQ_{\infty}) &= \int\sum\limits_{A\in\SP} -\psi_A(z)\log \psi_A(z)
d\mu(z) \\
&= \sum\limits_{A\in\SP}\int_{f(A)}\frac{1}{J_{\mu}f(y_A)}\log J_{\mu}f(y_A)
d\mu(z).
\end{split}
\end{equation*}
Since $f$ is injective on $A$, the definition of Jacobian implies
\begin{equation*}
H_{\mu}(\SP|\SQ_{\infty})=\sum\limits_{A\in\SP}\int_A \log J_{\mu}f(y) d\mu(y) =
\int\log J_{\mu}f(y) d\mu(y).
\end{equation*}
This concludes the proof.
\end{proof}

Now, it is an easy matter to complete the proof of the
proposition~\ref{p.rokhlin}. Indeed, since $\SP$ is generating, $h_{\mu}(f) =
h_{\mu}(f,\SP)$. On the other hand, it is well-known that $h_{\mu}(f,\SP) =
H_{\mu}(\SP|\SQ_{\infty})$ (see~\cite[Theorem 4.14]{W}). Therefore, the
lemma~\ref{l.rokhlin2} allow us to conclude $h_{\mu}(f)=\int\log J_{\mu}f d\mu$. 
\end{proof}

For later use, we show the following abstract lemma relating positive Lyapounov
exponents and hyperbolic times:

\begin{lemma}\label{l.appendixIII} Given any ergodic measure $\eta$ whose
Lyapounov exponents are all bigger than $8c$ then there exists some
$N\in\natural$ such that $f^N$ has infinitely many $c$-hyperbolic times
for $\eta$ almost every point.
\end{lemma}

\begin{proof}Since the Lyapounov exponents of $\eta$ are all bigger than
$8c$, for almost every $x\in M$, there exists $n_0(x)\geq 1$ such that
\begin{equation*}
\|Df^n(x)w\|\geq e^{6cn}\|w\|, \quad \text{ for all }w\in T_xM \text{ and }
n\geq n_0(x).
\end{equation*}
In other words, $\|Df^n(x)^{-1}\|\leq e^{-6cn}$ for every $n\geq n_0(x)$. Define
$\alpha_n=\eta (\{x:n_0(x)> n\})$. Because $f$ is a local diffeomorphism, we may
also fix $K>0$ such that $\|Df(x)^{-1}\|\leq K$ for all $x\in M$. Then,
\begin{equation*}
\int_M\log\|Df^n(x)\| d\eta\le -6cn+Kn\alpha_n = -(6c-K\alpha_n)n.
\end{equation*}
Since $\alpha_n\to 0$ when $n\to\infty$, by choosing $N$ large enough, we
guarantee that
\begin{equation*}
\int_M\frac{1}{N}\log\|Df^N(x)^{-1}\| d\eta <-4c<0.
\end{equation*}
Hence, the ergodicity of $\eta$ implies
\begin{equation*}
\lim\limits_{n\to\infty}\frac{1}{n}\sum\limits_{j=0}^{n-1}
\frac{1}{N}\log\|Df^N(f^{Nj}(x))^{-1}\| = 
\int_M\frac{1}{N}\log\|Df^N(x)^{-1}\| d\eta <-4c.
\end{equation*}
Now, the proof of the lemma is completed by taking $A=\sup\limits_{x\in M} -
\log\|Df(x)^{-1}\|$, $c_1=2c$, $c_2=3c$ and $b_j=-\log\|Df^N(f^{Nj}(x))^{-1}\|$
in lemma~\ref{l.Pliss}.
\end{proof}

\begin{corollary}\label{c.rokhlinA}
Under the hypothesis of theorem~\ref{t.A}, any measure $\eta$ with
$h_{\eta}(f)+\int\phi \; d\eta\geq\log\lambda$ satisfies Rokhlin's formula:
$$
h_{\eta} (f) = \int J_{\eta} f d\eta.
$$
\end{corollary}

\begin{proof} The proposition~\ref{p.rokhlin} says that it suffices to prove
that any $\eta$ with $h_{\eta}(f)+\int\phi d\eta\geq\log\lambda$ admits a 
$f$-generating partition. Without loss of generality, we can assume that $\eta$
is ergodic. 

The first step is to show that the Lyapounov exponents 
$\lambda_{1}(\eta)\geq\dots\geq\lambda_d(\eta)$ of $\eta$ are 
positive. 

We claim that $\lambda^d(\eta)\geq 
\log k-\max\limits_{x\in M}\log\|\Lambda^{d-1} Df(x)\|-\epsilon_0(f):=\kappa_0$ 
whenever
$h_{\eta}(f)+\int\phi d\eta\geq\log\lambda$. Indeed, this is an immediate
consequence of Ruelle's inequality:
\begin{equation*}
\log k\leq\log\lambda\leq h_{\eta} (f)+\int\phi d\eta\leq \lambda^d(\eta) + 
\max\limits_{x\in M}\log\|\Lambda^{d-1} Df(x)\|+\epsilon_0(f)
\end{equation*} 

Finally, the conclusion of the proof of the corollary~\ref{c.rokhlinA} is: we 
know that the Lyapounov exponents of any ergodic measure $\eta$ are all bigger 
than $\kappa_0$; so, the lemma~\ref{l.appendixIII} above means that
$\eta$-almost every point has infinitely many hyperbolic times with respect to
$f^N$, for some $N\in\natural$. 

Hence, we can use the lemma~\ref{l.3.5} to obtain that any partition $\SP$ with
diameter sufficiently small is generating. This completes the proof.
\end{proof}

\begin{corollary}\label{c.rokhlinB}
Under the hypothesis of theorem~\ref{t.B}, any equilibrium
measure $\eta$ satisfies Rokhlin's formula:
$$
h_{\eta} (f) = \int J_{\eta} f d\eta.
$$
\end{corollary}

\begin{proof}We give two proofs of this result. 

The first proof is based on the
positivity of all Lyapounov exponents of any ergodic equilibrium state $\eta$ 
(see theoerem~\ref{t.K}). So, the abstract lemma~\ref{l.appendixIII} can be used
to conclude the existence of generating partitions, and then, the Rokhlin's
formula is a consequence of the abstract proposition~\ref{p.rokhlin}.

The second proof consists in a transference of the problem of Rokhlin's formula 
for measures on $M$ to the same problem in some subshift of finite type, where 
this formula is known to be true.\footnote{In fact, although these two proofs
are formally different, they are the same in spirit since the proof of Rokhlin's
formula for invariant measures in a subshift of finite type goes in the lines of
lemma~\ref{l.rokhlin}.} Consider the partition $\SR$ and define 
$\Pi:M\to\Sigma^+$ a map over a subshift of finite type $\Sigma^+$ by 
$$\Pi (x) = (i_0,\dots,i_n,\dots) \quad \text{ such that } f^n(x)\in R_{i_n}.$$
Observe that this map is an ergodic conjugacy with respect to any measure $\eta$
such that $[x] = \{x\}$ at $\eta$-almost every $x$. Here, if $(i_n)$ is the
itinerary of $x$ (i.e., $f^n(x)\in R_{i_n}$), then   
\begin{equation*}
[x]:=[i_0,\dots,i_n,\dots]:=\{y\in M: f^n(y)\in R_{i_n}\}
\end{equation*}  
Note that any measure with infinitely many hyperbolic times at almost every
point satisfies $[x]=\{x\}$ at generic points. So, any equilibrium measure
$\eta$ can be transported to shift as $\Pi^*\eta$. Therefore, $\eta$ satisfies
Rokhlin's formula, since $\Pi^*\eta$ verifies it (see~\cite{BS} for more
details). 
\end{proof}

\bibliographystyle{alpha}
\bibliography{bib}

\vspace{1.5cm}

\noindent 	Alexander Arbieto ({\tt alexande{\@@}impa.br}) \\
		Carlos Matheus ({\tt matheus{\@@}impa.br}) \\
		IMPA, Estrada D. Castorina 110, Jardim Bot\^anico, 22460-320 \\ 
		Rio de Janeiro, RJ, Brazil \\

\end{document}